\documentclass[review,11pt]{elsarticle}

\journal{Fuzzy Sets and Systems}
\usepackage[utf8]{inputenc}
\usepackage[T1]{fontenc}
\usepackage{enumerate,enumitem}
\usepackage{amsfonts,amsmath,amssymb,amsthm,geometry,wrapfig,url,yfonts,fancyhdr}
\usepackage{subfigure,graphicx,color,mathrsfs }
\everymath{\displaystyle}
\usepackage{ulem}
\usepackage{dsfont}
\usepackage[colorlinks=true]{hyperref}

\newtheorem{tw}{Theorem}[section]	
\newtheorem{corollary}[tw]{Corollary}

\newtheorem{proposition}[tw]{Proposition}
\newtheorem{lemma}[tw]{Lemma}
\theoremstyle{definition} 
\newtheorem{definition}[tw]{Definition}
\newtheorem{example}[tw]{Example}
\newtheorem{remark}[tw]{Remark}
\newproof{pf}{Proof}
\newproof{pot1}{\textbf{Proof of Theorem~\ref{tw:well}}}
\newproof{pot2}{\textbf{Proof of Theorem~\ref{tw:mon_2}}}
\newproof{pot3}{\textbf{Proof of Theorem~\ref{tw:mon_3}}}
\newproof{pot4}{\textbf{Proof of Theorem~\ref{tw:mon_n}}}

\newcommand\mR{{\mathbb R}}

\newcommand\mN{{\mathbb N}}

\newcommand\bM{\mathbf{M}} 

\newcommand{\cC}{\mathcal{C}}
\newcommand{\cCa}{\cC^{\mathrm{I}}}
\newcommand{\cCb}{\mathcal{C}^{\mathrm{II}}}

\newcommand\ag{\mathrm{Ag}}

\newcommand\sP{\mathsf{P}}

\newcommand{\sK}{\mathsf{K}}

\renewcommand\int{\mathsf{Int}_{\mR_+}}
\newcommand\sLi{\mathsf{Int}_{[0,1]}}

\newcommand{\sU}{\mathsf{U}}

\renewcommand\ge{\geqslant}
\renewcommand\le{\leqslant}

\newcommand\les{\preceq_{\small{\mathrm{Int}}}}

\newcommand\leqa{\preceq}
\newcommand\lea{\prec}

\newcommand\leqA[1][n]{\preceq_{(#1)}}

\newcommand\leqp{\preceq_{\mathrm{p}}} 
\newcommand\lep{\prec_{\mathrm{p}}} 

\newcommand\lev{\preceq_{\mathrm{V}}}

\newcommand\leqab{\preceq_{(\alpha,\beta)}} 

\newcommand\bx{\mathbf{x}}
\newcommand\by{\mathbf{y}}
\newcommand\bz{\mathbf{z}}
\newcommand\bv{\mathbf{v}}
\newcommand\bu{\mathbf{u}}

\newcommand\bX{\mathbf{X}}

\newcommand\bZ{\mathbf{Z}}

\newcommand{\rd}{\mathrm{d}}
\newcommand{\rF}{\mathrm{F}}
\newcommand{\rG}{\mathrm{G}}
\newcommand{\rH}{\mathrm{H}} 
\newcommand{\rL}{\mathrm{L}}
\newcommand{\rP}{\mathrm{P}}

\newcommand{\bzer}{\mathbf{0}} 
\newcommand{\bjed}{\mathbf{1}} 

\newcommand{\tx}{\textstyle{}}

\usepackage{xcolor}\definecolor{darkgreen}{rgb}{0,0.5,0}

\begin{document}
\begin{frontmatter}
\title{
The Choquet-like operator with respect to an admissible order as a tool for aggregating multivalued data}

\author{Micha{\l} Boczek\corref{cor1}\fnref{label1}}
\ead{michal.boczek.1@p.lodz.pl}
\cortext[cor1]{michal.boczek.1@p.lodz.pl}

\author{Tomasz Józefiak\fnref{label1}}
 \ead{tomasz.jozefiak@dokt.p.lodz.pl}

\author{Marek Kaluszka\fnref{label1}}
 \ead{marek.kaluszka@p.lodz.pl}

\author{Andrzej Okolewski\fnref{label1}}
\ead{andrzej.okolewski@p.lodz.pl}

\address[label1]{Institute of Mathematics, Lodz University of Technology, 93-590 Lodz, Poland}

\begin{abstract}
In this paper, we propose a~new generalization of the classical discrete Choquet integral to the multivalued framework in terms of an admissible order that refines the natural partial order on the considered value set. 
The new Choquet-like operator takes as input a~finite number of values of a~given type, in particular   real numbers, intervals, and vectors, and returns a~single output value of the same type as the input values.
We give necessary and sufficient conditions for the operator to be monotone with respect to the admissible order. 
We then provide a~complete characterization of the Choquet-like operator as an aggregation function with respect to the admissible order and study its selected special cases.
\end{abstract}

\begin{keyword}
Admissible order\sep
Admissible permutation\sep
Aggregation function\sep
Choquet-like operator \sep
Monotone measure\sep
Partial order
\end{keyword}
\end{frontmatter}

\section{Introduction}






%

Information fusion is central to both theoretical and applied research. In recent years, there has been a~growing interest in the application of fuzzy integrals to data fusion, with a~particular focus on the discrete Choquet integral. The integral has proven to be extremely useful in various machine learning problems \cite{Barrenechea2013, Dias2018, Grabisch}. This has led to several generalizations of the Choquet integral in the real number setting being proposed in the literature, resulting in concepts such as 
the CQO-integrals  \cite{batista2022}, $\mathrm{C}_{\mathrm{T}}$-integrals \cite{dimuro2020}, 
CC-integrals  \cite{lucca2017}, 
d-Choquet integrals \cite{bustince2021}, 
and $\mathrm{C}_{\mathrm{F_1F_2}}$\,-\,integrals \cite{lucca2019}. 
A detailed survey of many of these integrals can be found in \cite{dimuro2020}. 
A~unified approach to integration in the Choquet sense in the real number setting was presented in \cite{BHK22,BK23},  introducing the concept of an operator based on conditional aggregation and relation operators.     
This operator 
includes many Choquet-like integrals known in the literature,  in particular those in \cite{batista2022, dimuro2020, lucca2017, bustince2021, lucca2019}.

Recent research has increasingly focused on the generalization of classical fuzzy set theory \cite{Zadeh} to include interval-valued fuzzy sets \cite{pekala}  
and $n$-dimensional fuzzy sets \cite{Bedregal2012}. The degrees of membership for such extensions are certain sets that consist of two or more real numbers. 
Several generalizations of the Choquet integral in multivalued settings have been proposed, 
in particular for the interval setting \cite{bustince2013b}, the $n$-dimensional vector setting \cite{ferrero23, JKO24}, and the $n$-dimensional interval setting \cite{milfont2021}. 
The construction of such extensions requires, among other things, overcoming problems arising from the fact that the natural product order on a~set of vectors with real coordinates is not a~linear order. 
Bustince et al. \cite{bustince2013a} tackled this problem in the interval-valued setting by introducing the concept of 
admissible orders, i.e., total orders that refine the product order.

The main reason for using the standard Choquet integral \cite{choquet1953} for data aggregation is its native properties. 
It is desirable that the functionals that are extensions of the Choquet integral inherit its original 
properties. 
Monotonicity, the fundamental axiom of the aggregation function, is one of these properties.
Its verification is a~major challenge in determining whether an operator qualifies as an aggregation function. 
Among the various generalized Choquet integrals proposed in the literature, only some satisfy the monotonicity condition, while some others fulfill a~weakened monotonicity condition, such as the directional monotonicity or ordered directional monotonicity \cite{dimuro2020}.
A~full characterization of the monotonicity condition for two important variants of the extended Choquet-like operator in the 
real number framework, which as special cases include the operators studied in \cite{bustince2020, horanska2018, mesiar2016}, 
along with the necessary and sufficient conditions that determine when the analyzed operators qualify as aggregation functions, are presented in \cite{BJKO23}.

In this paper we focus on the aggregation of a~finite number of data of a~given type, e.g., real-valued, interval-valued, or $n$-dimensional-valued, from the perspective of Choquet integral theory.
We introduce the notion of a~Choquet-like operator with respect to an admissible order, which refines the natural order on the considered data set \cite{bustince2013a, milfont2021}. 
This operator extends the discrete Choquet integral to the multivalued framework, replacing subtraction and multiplication by a~general function, and addition by an operation compatible with the admissible order.
We characterize the classes of all functions for which the Choquet-like operator  is: (a) well defined, i.e., invariant with respect to permutations of the input data, (b) monotone 
with respect to the admissible order. 
Finally, we give a~complete characterization of the Choquet-like operator as an aggregation function with respect to the admissible order, and study some special cases in detail.

The structure of the paper is as follows. In Section~\ref{sec:prel} we present some preliminary definitions and results which are necessary for the rest of the paper. 
In Section~\ref{sec:1} we introduce 
a~Choquet-like operator  with respect to an admissible
order and study its 
well-definedness and monotonicity. In Section~\ref{sec:2} we provide necessary and sufficient conditions for the operator to be an 
aggregation function, and analyze these conditions in detail for two important particular variants of the operator.
To improve readability, we move the well-definedness proof and all related monotonicity proofs to Appendix.

\section{Preliminaries}\label{sec:prel}

Let $\mR_+$ and $\mN$ be the sets of non-negative real numbers and positive integers, respectively. 
By $[n]$ we denote the set $\{1,\ldots,n\}$ for $n\in \mN$. 
Let $\sU$ be a~non-empty universal set,  
sometimes called a~reference set or  universe of discourse.

\subsection{ \textbf{Ordered sets}}

A~\textit{partial order} $\leqp$ on $\sU$ is a~binary relation on $\sU$, which is:
\begin{itemize}[noitemsep]
    \item reflexive: $\bx\leqp \bx$ for any $\bx\in \sU$, 
    
    \item antisymmetric: for all $\bx,\bz\in\sU$, $\bx\leqp \bz$ and $\bz\leqp \bx$ implies $\bx=\bz$,

    \item transitive: for all $\bx,\by,\bz\in \sU$, if $\bx\leqp \by$ and $\by\leqp \bz$,  then $\bx\leqp \bz$.
\end{itemize}

In the rest  of the paper, we will write ``partial order'' instead of ``partial order on $\sU$'' when the form of the universal set $\sU$  is clear from the context.

Using the partial order $\leqp$,  we establish the relation $\lep$ as follows: for any elements $\bx$ and $\bz$ of the set $\sU$, the statement $\bx\lep \bz$ holds true if and only if $\bx\leqp \bz$ and $\bx\neq \bz$.
By the reflexivity and antisymmetry of the partial order $\leqp$, it follows that for any $\bx,\bz\in \sU$, $\bx = \bz$ if and only if $\bx \leqp \bz$ and $\bz \leqp \bx$. 
Hence, for all $\bx,\bz\in \sU$, $\bx \leqp \bz$ if and only if $\bx \lep \bz$ or $\bx =\bz$.

A~non-empty set $\sK\subseteq \sU$ is said to be \textit{bounded w.r.t.~a~partial order} $\leqp$ if it has both the least element $\bzer$ and the greatest
element $\bjed$ w.r.t.~$\leqp$, that is, 
$\bzer,\bjed \in \sK$ and $\bzer\leqp \bx \leqp \bjed $ for all $\bx \in \sK$.
A~partial order $\leqa$ for which every  two elements are comparable, i.e., 
$\bx\leqa \bz$ or $\bz\leqa \bx$ for any $\bx,\bz\in\sU$, is called a~\textit{total order}.  


We will now present some examples of restricted subsets for certain universal sets with usual partial orders, which  will be used later in the paper.

\begin{example}\label{ex:1}
\begin{enumerate}[label=(\alph*), noitemsep]
    \item For $\sU=\mR_+$, 
     the set $\sK=[0,1]$ is bounded w.r.t.~the usual order~$\le$ with the least element $0$ and the greatest element $1.$ 
     Hereafter, the symbol $[a, b]$ with $a\le b$ represents an interval in $\mR_+$. 
         
    \item Let $\sU = \int$ be the family of all closed subintervals of $\mR_+$. 
    The natural partial order is the \textit{interval order} $\les$ defined as follows
    \begin{align*}
    \bx \les \bz\qquad \Leftrightarrow\qquad x^l\le z^l \text{ and } x^u\le z^u,  
    \end{align*}
    for any $\bx=[x^l, x^u]\in \int$ and $\bz=[z^l, z^u]\in \int$.
    Then  $\sK=\sLi$, where $\sLi=\{\bx = [x^l, x^u]\mid 0\le x^l\le x^u\le 1\}$, is a~bounded set w.r.t.~$\les$ with the least element $[0,0]$ and the greatest element $[1,1]$.
        
    \item Let $\sU = \mR^k_+$ with $k\in\mN$. 
    The usual partial order $\lev$ is defined as follows  
    \begin{align*}
        \bx \lev \bz \quad \Leftrightarrow\quad x_i \le z_i \text{ for all } i\in [k],
    \end{align*}
    for any $\bx = (x_1,\ldots,x_k)^T\in \mR^k_+$ and $\bz=(z_1,\ldots,z_k)^T\in \mR^k_+$, where $\bx^T$ denotes the transpose of a~vector  $\bx$.
    Then $\sK= [0,1]^k$ is a~bounded set w.r.t.~$\lev$ with the least element $(0,\ldots,0)^T$ and the greatest element $(1,\ldots,1)^T$. 
\end{enumerate}
\end{example}

\begin{definition}
Let $\leqp$ be a~partial order.
A total order $\leqa$ is called an \textit{admissible order  w.r.t.~$\leqp$} if for any $\bx,\bz\in \sU$, it holds $\bx \leqa \bz$ whenever $\bx \leqp \bz$. 
\end{definition}

The notion of the admissible order was introduced in \cite{bustince2013a} in the context of comparing two intervals. 
One of the most significant admissible orders  w.r.t.~$\les$ is the $(\alpha, \beta)$-order  defined as follows:
for any $\bx,\bz \in \int$
\begin{align*}
\bx \preceq_{(\alpha, \beta)} \bz \quad \Leftrightarrow \quad 
\begin{cases}
    K_{\alpha}(\bx) < K_{\alpha}(\bz) \text{ or } \\K_{\alpha}(\bx) = K_{\alpha}(\bz) \text{ and } K_{\beta}(\bx) \le K_{\beta}(\bz),
\end{cases}
\end{align*}
where $\alpha, \beta \in [0, 1]$, $\alpha \neq \beta$, and $K_{\alpha}(\by) = (1-\alpha) y^l + \alpha y^u$ for any $\alpha \in [0,1]$ and $\by\in \int$.
The $(\alpha, \beta)$-order generalizes three well-known admissible orders:
(i) the lexicographical order for $(\alpha, \beta) = (0,1)$,
(ii) the antilexicographical order for $(\alpha, \beta) = (1,0)$,
(iii) the Xu-Yager order for $(\alpha, \beta) = (0.5, 1)$. 
In \cite{BJK21} it is proved that in a~certain class of orders, the only admissible order  w.r.t.~$\les$ is the $(\alpha, \beta)$-order.



\subsection{\textbf{Addition operation on $\sU$}}

A~binary operation $\oplus\colon \sU\times \sU \to \sU$ is called {\it addition} if it satisfies the following properties:
\begin{enumerate}[noitemsep, label=(\alph*)]
    \item commutativity: $\bx \oplus \bz = \bz \oplus \bx$ for any $\bx,\bz$, 
    \item associativity: $(\bx \oplus \by) \oplus \bz = \bx \oplus (\by \oplus \bz)$ for any $\bx,\by,\bz$. 
\end{enumerate}

Note that $(\sU, \oplus)$ is a~commutative semigroup.  In particular, $([0,1], \mathrm{T})$ and $([0,1], \mathrm{S})$ are commutative semigroups for any $t$-norm $\mathrm{T}$ and $t$-conorm $\mathrm{S}$ \cite{fodor2000}.
For more examples of semigroups, we refer to \cite{robinson2003}. Examples of addition operations $\oplus$ for the special cases of $\sU$ described in Example~\ref{ex:1} are given in Table~\ref{tab:1}.

\begin{table}[h!]
\centering
\begin{tabular}{|c| c| c|}\hline
 $\sU$ & $\oplus$ & $\mathbf{x \oplus z}$ \\ \hline
$\mR_+$ & $+$ & $x+z$ (standard addition) \\ \hline
 $\int$  & $\oplus_{\mathrm{iv}}$ & $[x^l,x^u]\oplus_{\mathrm{iv}} [z^l,z^u]=[x^l+z^l, x^u + z^u]$ \\ \hline
$\mR^k_+$  & $\oplus_{\mathrm{vv}}$ & $(x_1,\ldots,x_k)^T\oplus_{\mathrm{vv}}(z_1,\ldots,z_k)^T=(x_1+z_1,\ldots, x_k+z_k)^T$  \\ \hline
\end{tabular}
    \caption{Examples of addition  operations  for selected cases of $\sU$.}
    \label{tab:1}
\end{table}


\begin{definition}\label{def:2.4}
A partial order $\leqp$  
is said to be \textit{compatible with $\oplus$} if for any $\bx_1,\bx_2,\bv\in \sU$  we have 
$$
\bx_1 \leqp \bx_2\qquad \Rightarrow\qquad \bx_1 \oplus \bv \leqp \bx_2 \oplus \bv.
$$
Replacing $\leqp$ with $\lep$, we get the notion of \textit{compatibility of $\lep$ with $\oplus$}. 
\end{definition}


All  partial orders from Example~\ref{ex:1} are compatible with the corresponding operations $\oplus$ from Table~\ref{tab:1}.

\begin{definition}
 An addition operation $\oplus$  has the \textit{cancellation law}, whenever 
for any $\bx_1,\bx_2,\bv \in \sU$  the following implication holds
\begin{align*}
\bx_1 \oplus \bv = \bx_2 \oplus\bv \qquad \Rightarrow \qquad\bx_1 = \bx_2.
\end{align*}
\end{definition}

All operations from Table~\ref{tab:1} have the cancellation law. 
We now present some auxiliary results.

\begin{lemma}\label{lem:1}
Let $\leqp$ be a~partial order.
If  $\lep$  is compatible with $\oplus$, then $\leqp$ is compatible with $\oplus$.
\end{lemma}
\begin{pf}
By assumption, we have $\bx_1 \oplus \bv \leqp \bx_2 \oplus \bv$ for any $\bx_1\lep \bx_2$ and $\bv\in \sU$.  Clearly, if $\bx_1 = \bx_2$, then $\bx_1 \oplus \bv \leqp \bx_2 \oplus \bv$ as $\bx_1 \oplus \bv=\bx_2 \oplus \bv$.   \qed
\end{pf}

The minimum operation  on $\mR_+$ is commutative and  associative, and  the usual order $\le$ is compatible with this operation.
However, the relation $<$ is not compatible with the minimum operation as $\min\{x_1, c\} = \min\{x_2, c\}$ for any $0 \le c < x_1 < x_2$.

To reverse the implication of Lemma~\ref{lem:1}, an additional assumption is necessary.

\begin{lemma}\label{lem:2}
Let $\oplus$ satisfy the cancellation law and $\leqp$ be a~partial order. 
If $\leqp$  is compatible with $\oplus$, then $\lep$ is compatible with $\oplus$.
\end{lemma}
\begin{pf}
By  assumption, $\bx_1 \oplus \bv \leqp \bx_2 \oplus \bv$ for any $\bx_1 \lep \bx_2$ and $\bv$. 
Suppose there exists $\bx_1 \lep \bx_2$ and $\bv$ such that $\bx_1 \oplus \bv = \bx_2 \oplus \bv$. Then by the cancellation law, we have $\bx_1= \bx_2$, a~contradiction.  
The proof is complete. 
\qed
\end{pf}

\begin{lemma}\label{lem:3}
Assume that $\leqa$ is  an admissible order w.r.t.~$\leqp$ such that  $\lea$  is compatible with $\oplus$.
Then: 
\begin{enumerate}[noitemsep, label=(\alph*)]
    \item  $\oplus$ has the cancellation law, 

    \item  for any $\bx_1,\bx_2,\bv\in\sU$, if $\bx_1 \oplus \bv \leqa \bx_2 \oplus \bv$,     then $\bx_1 \leqa \bx_2$,

    \item  for any  $\bu,\bv\in \sK$, if 
     $\bzer$ is the least element of $\sK$ w.r.t.~$\leqp$,
    $\bzer \oplus \bzer \in \sK$ and $\bu \oplus \bv=\bzer$, then  $\bu=\bzer$ and $\bv=\bzer$. 
\end{enumerate}
\end{lemma}
\begin{pf}
(a) Suppose that $\bx_1 \oplus \bv = \bx_2 \oplus \bv$ for  some $\bv$  and $\bx_1 \neq \bx_2$. 
Since $\leqa$ is a~total order, any two elements are comparable, i.e.,  for any $\bx,\by\in \sU$ we have either $\bx \lea\by$ or $\bx=\by$ or $\by\lea \bx$. 
In our case, we get either $\bx_1 \lea \bx_2$ or $\bx_2 \lea \bx_1$, as $\bx_1 \neq \bx_2$. 
If $\bx_1 \lea \bx_2$, then by the compatibility of $\lea$ we obtain 
$\bx_1 \oplus \bv\lea \bx_2 \oplus \bv$, so it is not true that $\bx_1 \oplus \bv = \bx_2 \oplus \bv$, a~contradiction. 
The case $\bx_2 \lea \bx_1$ can be examined in a~similar way. 

(b)  
Suppose that $\bx_1 \oplus \bv \leqa \bx_2 \oplus \bv$  and $\bx_2 \lea \bx_1$ for  some  $\bx_1,\bx_2,\bv$. By  the compatibility of $\lea$, we get 
$\bx_2 \oplus \bv\lea \bx_1 \oplus \bv$, a~contradiction.

(c)
Since $\bu \oplus \bv =\bzer$, $\bzer \oplus \bzer \in \sK$, and $\bzer$ is the least element of $\sK$, we have $\bu \oplus \bv =\bzer \leqa  \bzer \oplus \bzer.$ 
From Lemma~\ref{lem:1}, we get $\bzer \oplus \bzer \leqa \bu \oplus \bzer.$ 
Thus, by the transitivity of $\leqa$, $\bu \oplus \bv \leqa \bu \oplus \bzer.$ 
Point (b) now gives  $\bv \leqa \bzer.$ Hence, $\bv =\bzer$.
Similarly, it can be shown that $\bu = \bzer.$ 
\qed
\end{pf}

\subsection{\textbf{$\sK$-valued aggregation function}}

Let $\leqa$ be an admissible order  w.r.t.~$\leqp$.
Fix $\bX = (\bx_1,\ldots,\bx_n)\in \sU^n$ and $\bZ= (\bz_1,\ldots, \bz_n)\in \sU^n$. 
We define a~binary relation on $\sU^n$, say $\leqA$, as 
\begin{align*}
    \bX \leqA \bZ\qquad \Leftrightarrow \qquad \bx_i\leqa \bz_i \text{ for any } i\in [n].
\end{align*}




\begin{definition}\label{def:ag}
Let $n\ge 2$, $\leqa$ be an admissible order  w.r.t.~$\leqp$, and $\sK$ be a~bounded set w.r.t.~$\leqp$ with the least element $\bzer$ and the greatest element $\bjed$.
A~mapping $\ag\colon \sK^n \to \sK$ is called an (\textit{$n$-ary}) \textit{$\sK$-valued aggregation function w.r.t.~$\leqa$} if it satisfies the following monotonicity condition:
\begin{align*}
   \ag(\bX)\leqa \ag(\bZ)\qquad \text{if} \qquad \bX\leqA \bZ,
\end{align*}
and meets the boundary conditions:
\begin{align}\label{ag:2}
    \ag(\bzer,\ldots, \bzer)=\bzer\quad \text{and} \quad \ag(\bjed,\ldots, \bjed)=\bjed.
\end{align}
\end{definition}

Putting $\sK = [0,1]$ and $\le$ in place of $\leqa$
in Definition~\ref{def:ag}, we get the  well-known definition of the real-valued 
aggregation function
(see \cite{calvo2002}). 
If we consider 
the set $\sK = \sLi$ with the partial order $\les$ we recover the concept of an~interval-valued (IV, for short) aggregation function (see \cite{bentkowska2020}). 
In the case of $\sK=[0,1]^k$, to 
the best of our knowledge, the definition of the vector-valued (VV, for short) aggregation function as a~special form of Definition~\ref{def:ag} has not yet been considered in the literature (although a~definition of the monotonicity w.r.t.~the admissible order $\leqa$ 
has been introduced in \cite[Def.~4.5]{sara2023}). 
On the other hand, in the case of the $k$-dimensional upper simplex $\sK = \{(x_1,\ldots, x_k)\mid  0\le x_1 \le x_2\le \ldots \le x_k\le 1\}$ such a~definition has been implemented in the context of $k$-dimensional fuzzy sets (see \cite[Def.~7]{milfont2021}). 

\subsection{\textbf{$\sK$-valued dissimilarity function}}\label{sec:dis}

\begin{definition}\label{def:dis}
Let $\leqa$ be an admissible order w.r.t.~$\leqp$  and $\sK$ be a~bounded set w.r.t.~$\leqp$ with the least element $\bzer$ and the greatest element $\bjed$.
A~function $\rd\colon \sK^2 \to \sK$ is said to be a~\textit{$\sK$-valued dissimilarity function w.r.t.~$\leqa$} if for all $\bx,\by,\bz \in \sK$,  it satisfies the following conditions:
\begin{enumerate}[noitemsep, label=(\alph*)]
    \item $\rd(\bx, \bz)=\rd(\bz, \bx)$,
    \item $\rd(\bzer, \bjed) = \bjed$,
    \item $\rd(\bx,\bx)= \bzer$,
    \item if $\bx \leqa \by$ and $\by\leqa \bz$, then $\rd(\bx,\by)\leqa \rd(\bx, \bz)$ and $\rd(\by, \bz) \leqa \rd(\bx,\bz)$.
\end{enumerate}
\end{definition}

For $\sK = [0,1]$ we get the definition of the \textit{real-valued dissimilarity function} proposed in ~\cite{bustince2008}, while for $\sK = \sLi$ we recover the definition of the \textit{IV dissimilarity function} introduced in ~\cite{takac2022}.
These functions are designed to quantify the dissimilarity between two input values in various settings where the difference causes problems \cite{bustince2021}.

\section{Main results}
\label{sec:1}

From now on, we assume that
a~set $\sK\subseteq \sU$ is bounded  w.r.t.~a~partial order $\leqp$ with the least element $\bzer$ and the greatest element $\bjed$ such that $\bzer \lep \bjed$ and that $\leqa$ is an admissible order  w.r.t.~$\leqp$.

By $\bM$ we denote the set of all \textit{capacities} on $[n]$, i.e., the set functions $\mu\colon 2^{[n]}\to [0,1]$  
satisfying $\mu(C)\leqslant\mu(D)$ whenever $C\subset D\subseteq [n]$ with $\mu(\varnothing)=0$ and $\mu([n])=1$. 
A~permutation $\sigma$ of $[n]$ is called \textit{admissible w.r.t.~$\leqa$} for $\bX = (\bx_1,\ldots,\bx_n)\in \sK^n$ if  $\bx_{\sigma(1)}\leqa \ldots\leqa \bx_{\sigma(n)}$. 
For a~fixed $\bX\in \sK^n$, by $\Pi_{\bX}$ we  denote the set of all admissible permutations w.r.t.~$\leqa$. 
Clearly, $\Pi_{\bX} \neq \varnothing$.

\begin{definition}\label{def:choquet}
Let $\rL\colon \sK^2 \times [0,1]^2 \to \sK$.
The \textit{Choquet-like operator w.r.t.~$\leqa$}
(\textit{$\cC_{\leqa}$\,-\,operator}, in short)  of  $\bX=(\bx_1,\ldots,\bx_n)\in \sK^n$  and~$\mu\in \bM$ is defined as follows
\begin{align}\label{dc}
    \cC_{\leqa, \sigma, \rL}(\bX, \mu) =\bigoplus_{i=1}^n \rL(\bx_{\sigma(i)}, \bx_{\sigma(i-1)}, \mu(B_{\sigma(i)}),\mu(B_{\sigma(i+1)})),
\end{align}
where  $\sigma\in \Pi_{\bX}$,  $\bx_{\sigma(0)}=\bzer$, $B_{\sigma(i)}= \{\sigma(i),\ldots,\sigma(n)\}$ for $i\in [n]$, and $B_{\sigma(n+1)}=\varnothing$.
Hereafter, ${\tx \bigoplus_{i=1}^n \bz_i = \bz_1 \oplus \ldots \oplus \bz_n.}$ 
\end{definition}

For $\sU=\mR_+$, $\sK=[0,1]$,  and $\oplus = +$, the $\cC_{\le}$\,-\,operator   is called the fuzzy extended Choquet-like operator in \cite{BK23}.

\begin{remark}
For a~fixed $n\in \mN$, if $\leqp$ is compatible with $\oplus$, 
we have
\begin{align*}
   \underbrace{\bzer \oplus\ldots \oplus \bzer}_n \leqp \cC_{\leqa, \sigma, \rL}(\bX, \mu)\leqp \underbrace{\bjed \oplus\ldots \oplus \bjed}_n
\end{align*}
for any $\bX,\mu,\sigma,\rL$.
\end{remark}

Notable  examples of the $\cC_{\leqa}$\,-\,operator, as explored in the literature \cite{BK23, bustince2013b, lucca2017, takac2022}, arise for:
\begin{enumerate}[noitemsep, label = (\alph*)]
    \item $\rL(\bx_1, \bx_2, b_1, b_2 ) = \mathrm{G}(\bx_1, b_1,b_2)$, 

    \item $\rL(\bx_1, \bx_2, b_1, b_2 ) = \mathrm{G}(\bx_1, \bx_2,b_1)$.
\end{enumerate}
In these cases we obtain $\cCa_{\leqa}$\,-\,operator and $\cCb_{\leqa}$\,-\,operator, respectively,
\begin{align*}
    \cCa_{\leqa, \sigma, \mathrm{G}}(\mathbf{X},\mu)& = \bigoplus_{i=1}^n \mathrm{G}(\bx_{\sigma(i)}, \mu(B_{\sigma(i)}),\mu(B_{\sigma(i+1)})), \\
    \cCb_{\leqa, \sigma, \mathrm{G}}(\mathbf{X},\mu)& = \bigoplus_{i=1}^n \mathrm{G}(\bx_{\sigma(i)}, \bx_{\sigma(i-1)}, \mu(B_{\sigma(i)})).
\end{align*}

For the sake of brevity, when considering $\sK = [0,1]$, $\sK = \sLi$, or $\sK = [0,1]^k$, we will assume the natural partial orders given in Example~\ref{ex:1}, as well as $\sU$ with  $\oplus$ given in Table~\ref{tab:1}. 
For example, if $\sK = \sLi$, then  $\sU= \int$, $\leqp \,= \,\les$, and $\oplus = \oplus_{\mathrm{iv}}$.

\begin{remark}\label{rem:3.3}
We will now give a~brief overview of what has been done for  the
$\cCa_{\leqa}$\,-\,operator and the $\cCb_{\leqa}$\,-\,operator 
from the perspective of  $\sK$-valued aggregation functions, along  with some applications.
\begin{enumerate}[noitemsep, label = (\alph*)]
    \item 
    Consider $\sK= [0,1]$. 
    The properties of the $\cCa_{\leqa}$\,- and $\cCb_{\leqa}$\,-\,operators, such as well-definedness and monotonicity, have been investigated in \cite{BJKO23, BK23}. In  \cite{BJKO23}, the authors characterized all functions $\rG$ for which these operators qualify as real-valued aggregation functions. 
    In particular, they provided characterizations for certain cases: (i) for the $\cCa_{\leqa}$\,-\,operator with $\rG(a, b_1, b_2) = \mathrm{F}(a, b_1 - b_2)$ for $b_1\ge b_2$ (as examined in \cite{horanska2018}), and (ii) for the $\cCb_{\leqa}$\,-\,operator with $\rG(a_1, a_2, b) = \rF(a_1 - a_2, b)$ for $a_1\ge a_2$ (as studied in \cite{lucca2016b, mesiar2016}). 
    Characterizations have also been developed for other special cases of the $\cCb_{\leqa}$\,-\,operator, including the CQO-integral \cite{batista2022}, 
    $n$-ary discrete Choquet integral \cite{bustince2021}, 
    $\mathrm{C}_{\mathrm{T}}$-integral \cite{dimuro2020},
    CC-integral \cite{lucca2017}, 
    and CO-integral \cite{lucca2016}. 
    These operators can be used  for example in:
    multi-layer classifier ensembles \cite{batista2022},
    image reduction \cite{dias2019},
    multimodal fuzzy fusion decision \cite{ko2019},
    fuzzy rule-based classification systems \cite{lucca2016, lucca2016b}, and
    the classification of reflectance measurements from hyperspectral images of grapes \cite{lucca2018a}.

    \item For $\rL(\bx_1, \bx_2, b_1, b_2) = \rP(\rd(\bx_1,\bx_2), b_1)$  with $\sP\colon \sK \times [0,1] \to \sK$, $\rd\colon \sK^2 \to \sK$ such that $\rP(\bzer, c)=\bzer$ and $\rd(\bx,\bx)=\bzer$ for any $\bx,c$, we get a~special case of the $n$-ary discrete $d_G$-Choquet integral defined in \cite[Def.~3.1]{takac2022}. 
    For $\sK=\sLi$ and $\rd$ being  an IV dissimilarity function,
    the authors in \cite{takac2022} 
    provided  the condition  under which the $n$-ary discrete $d_G$-Choquet integral qualifies as an $n$-ary  IV aggregation function 
    w.r.t.~$\leqab$, cf. Section~\ref{Sec:4.2.1}.
    This approach is applied to aggregate the decisions made by multiple classifiers.

    \item  
    Let $\sK=\sLi$ and  $\rG(\bx, b_1, b_2) = (b_1- b_2)_+ \odot_{\mathrm{iv}} \bx$,  where  $\odot_{\mathrm{iv}}$ is defined in Table~\ref{tab:2} (see Section~\ref{sec:Gdis}) and $b_+ = \max\{b,0\}.$
    The  $\cCa_{\preceq_{A,B}}$\,-\,operator is the $\preceq_{A,B}$\,-\,Choquet integral \cite{bustince2013b,paternain2019}, where the admissible order $\preceq_{A,B}$  w.r.t.~$\les$ on $\sLi$ is defined in \cite[Prop.~3.2]{bustince2013a}.
    In \cite[Ex.~3.11]{BJK23} it was shown, that the  $\cCa_{\leqa}$\,-\,operator with an admissible order $\leqa$ w.r.t.~$\les$ is the special case of   $\mathrm{IVCSO}_{\leqa}$.

   \item If we assume that $\sK= [0,1]^k$ with $k \ge 2$ and $\rL(\bx_1, \bx_2, b_1, b_2) = (b_1 - b_2)_+ \odot_{\mathrm{vv}} \bx_1$, where  $\odot_{\mathrm{vv}}$ is defined in Table~\ref{tab:2} (see Section~\ref{sec:Gdis}),  we get a~vector-valued Choquet-like operator that has not been considered in the literature. 
   Another type of Choquet-like  operator for vector-valued data has been studied in \cite{ferrero23}.
\end{enumerate}
\end{remark}


We will now characterize the conditions under which the $\cC_{\leqa}$\,-\,operator is well defined. This property will be crucial in determining the conditions for the monotonicity of the operator. Applying the above results, in the next section we will give necessary and sufficient conditions for the  $\cC_{\leqa}$\,-\,operator to be a~$\sK$-valued aggregation function w.r.t.~$\leqa$.

\subsection{\textbf{Well-definedness}}\label{sec:wd}

The $\cC_{\leqa}$\,-\,operator is said  to be {\it well defined}  if, for  any fixed $\bX\in \sK^n$ and  $\mu\in \bM$  formula \eqref{dc} yields the same value for all  $\sigma\in\Pi_{\bX}$. 

\begin{tw}\label{tw:well}
Assume that $\oplus$ satisfies the cancellation law. 
\begin{itemize}
    \item The $\cC_{\leqa}$\,-\,operator with $n=2$ is well defined if and only if 
    
    \begin{enumerate}[label=(WD2), leftmargin=1.4cm]
    \item \label{WD2}
    $[0,1] \ni c \mapsto \rL(\bx, \bzer, 1, c) \oplus \rL(\bx, \bx, c, 0)$ is a~constant function  for any $\bx\in \sK$.
    \end{enumerate}

    \item The $\cC_{\leqa}$\,-\,operator with $n=3$ is well defined if and only if 
    
    \begin{enumerate}[label=(WD3), leftmargin=1.4cm]
    \item \label{WD3} 
    $[b, 1]\ni c\mapsto \rL(\bx, \bzer, 1, c) \oplus \rL(\bx, \bx, c, b)$ 
    and 
    $[0, b] \ni c\mapsto \rL(\bx_1, \bx_2, b, c) \oplus \rL(\bx_1, \bx_1, c, 0)$ 
    are  constant functions for any $\bx,\bx_i\in \sK$, $\bx_2 \leqa \bx_1$, and $b\in [0,1]$.
    
    \end{enumerate}

    \item The $\cC_{\leqa}$\,-\,operator with $n\ge 4$ is well defined if and only if 
    \begin{enumerate}[label=(WDn), leftmargin=1.4cm]
    \item \label{WDn}
    $[b_2, b_1] \ni c\mapsto \rL(\bx_1, \bx_2, b_1, c) \oplus \rL(\bx_1, \bx_1, c, b_2)$ is a~constant function for any $\bx_2\leqa \bx_1$ and $b_2\le b_1$.
    \end{enumerate}
\end{itemize}
\end{tw}
The proof of Theorem~\ref{tw:well} is postponed  to Appendix~A.

\subsection{\textbf{Monotonicity}}\label{sec:mon}

In this part we provide a~complete characterization  of the following monotonicity condition for the $\cC_{\leqa}$\,-\,operator 
\begin{enumerate}[label=(M), leftmargin=1.4cm]
    \item \label{m1}
    if $\bX,\bZ\in \sK^n$ and $\bX \leqA \bZ$, then 
    $\cC_{\leqa, \sigma, \rL}(\bX, \mu)\leqa \cC_{\leqa,\tau,\rL}(\bZ, \mu)$ 
    for any $\sigma\in \Pi_{\bX}$, $\tau\in \Pi_{\bZ}$, and $\mu\in\bM$.
\end{enumerate}



\begin{definition}
We say that a~function $\rH\colon \sK_0\to \sU$, where $\sK_0\subseteq \sK$, is \textit{non-decreasing w.r.t.~$\leqa$}  if $\rH(\bx_1) \leqa \rH(\bx_2)$ for $\bx_1\leqa \bx_2.$ 
\end{definition}

Let $[\bu, \bv] = \{\bx \in \sU \colon \bu \leqa \bx \leqa \bv\}$  with $\bu,\bv\in \sU$.
Denote by $[\bu, \bv]_{\sK} = [\bu, \bv] \cap \sK$.

We distinguish three cases: $n=2$, $n=3$, and $n\ge 4$.
Proofs of the corresponding theorems are given in Appendix B.

\begin{tw}\label{tw:mon_2}
Let $n=2$.
Assume that $\lea$ is compatible with $\oplus$.
The monotonicity condition~\ref{m1} holds if and only if the following conditions are satisfied:
\begin{enumerate}[noitemsep, label=(\alph*)]
    \item the $\cC_{\leqa}$\,-\,operator is well defined,
    
    \item 
    $[\bzer, \bv]_{\sK}\ni \bx \mapsto \rL(\bx, \bzer, 1, b) \oplus \rL(\bv, \bx,  b, 0 )$ is a~non-decreasing  function w.r.t.~$\leqa$ for any $\bv\in \sK$ and $b\in [0,1]$,

    \item $[\bu, \bjed]_{\sK} \ni\bx\mapsto \rL(\bx, \bu, b, 0)$ is a~non-decreasing function w.r.t.~$\leqa$ for any $\bu\in \sK$ and $b\in [0,1]$.
\end{enumerate}
\end{tw}

\begin{tw}\label{tw:mon_3}
Let $n=3$.
Assume that $\lea$ is compatible with $\oplus$. 
The monotonicity condition~\ref{m1}  is true if and only if
\begin{enumerate}[label=(\alph*), noitemsep]
    \item the $\cC_{\leqa}$\,-\,operator is well defined, 

    \item 
    $[\bzer, \bv]_{\sK}  \ni \bx\mapsto \rL(\bx, \bzer, 1, b_1) \oplus \rL(\bv, \bx, b_1, b_2)$ 
    is a~non-decreasing  function w.r.t.~$\leqa$ for any $\bv\in \sK$ and $b_2\le b_1$,

    \item 
    $[\bu, \bv]_{\sK}\ni \bx\mapsto \rL(\bx, \bu, b_1, b_2) \oplus \rL(\bv, \bx,  b_2, 0)$
    is a~non-decreasing  function w.r.t.~$\leqa$ for any $\bu\leqa \bv$ and $b_2\le b_1$,

     \item  
    $[\bu, \bjed]_{\sK}\ni\bx\mapsto \rL(\bx, \bu, b, 0)$ 
     is a~non-decreasing function w.r.t.~$\leqa$ for any $\bu\in \sK$ and $b\in [0,1]$.
\end{enumerate}
\end{tw}

\begin{tw}\label{tw:mon_n}
Let $n\ge 4$.
Assume that $\lea$ is compatible with $\oplus$. 
The monotonicity condition~\ref{m1}  is true if and only if 
\begin{enumerate}[label=(\alph*), noitemsep]
    \item the $\cC_{\leqa}$\,-\,operator is well defined, 

    \item $[\bu, \bv]_{\sK}\ni \bx\mapsto \rL(\bx, \bu, b_1, b_2) \oplus \rL(\bv, \bx,  b_2, b_3)$ 
    is a~non-decreasing function w.r.t.~$\leqa$ for any  $\bu\leqa \bv$ and $b_3\le b_2\le b_1$, 
    
    \item 
    $[\bu, \bjed]_{\sK}
    \ni\bx\mapsto \rL(\bx, \bu, b, 0)$ 
    is a~non-decreasing function w.r.t.~$\leqa$ for any $\bu\in \sK$ and $b\in [0,1]$.
\end{enumerate}
\end{tw}

\begin{remark}\label{rem:3.4}
By Lemma~\ref{lem:3}\,(a), condition (a) of Theorems~\ref{tw:mon_2}, \ref{tw:mon_3}, and \ref{tw:mon_n}, respectively, is equivalent to \ref{WD2}, \ref{WD3}, and \ref{WDn}.
\end{remark}

\section{Relationship to the $\sK$-valued aggregation function}
\label{sec:2}

According to Definition~\ref{def:ag}, the $\cC_{\leqa}$\,-\,operator is a~$\sK$-valued aggregation function if it satisfies the monotonicity condition~\ref{m1}, which is characterized in Theorems~\ref{tw:mon_2}\,-\,\ref{tw:mon_n}, and the boundary conditions \eqref{ag:2}. 

\begin{tw}\label{tw:ag}
Let $n\in \{2,3,\ldots\}.$
Assume  that $\lea$ is compatible with $\oplus.$
The $\cC_{\leqa}$\,-\,operator is an~$n$-ary $\sK$-valued aggregation function 
w.r.t.~$\leqa$ if and only if the following conditions hold:
    \begin{enumerate}[label=(\alph*), noitemsep]
            \item
             for $n=2$, $n=3$, and $n\ge 4$, respectively, the conditions of Theorems~\ref{tw:mon_2}, \ref{tw:mon_3}, and \ref{tw:mon_n} characterizing the monotonicity of the $\cC_{\leqa}$\,-\,operator are satisfied,

            \item  ${\tx\bigoplus_{i=1}^n \rL(\bzer,\bzer,b_i,b_{i+1}) = \bzer}$ for any $0=b_{n+1} \le b_n \le \ldots \le b_2 \le b_1=1$, 
            
            \item ${\tx \rL(\bjed,\bzer,b_1,b_2)\oplus \bigoplus_{i=2}^n \rL(\bjed,\bjed,b_i,b_{i+1})=\bjed}$  for any  $0=b_{n+1} \le b_n \le \ldots \le b_2 \le b_1=1$.       
    \end{enumerate}
\end{tw}

\begin{remark}\label{rem:42}
If we additionally assume that $\bzer \oplus \bzer\in \sK$,   then it follows from Lemma~\ref{lem:3}\,(c) that condition (b)  of Theorem~\ref{tw:ag} can be simplified to the condition: 
$\rL(\bzer,\bzer,b_1,b_{2}) = \bzer$ for any $b_2\le b_1.$
\end{remark}

In Sections \ref{sec:agC1} and \ref{sec:agC2}, 
we will present specific results regarding two significant operators: the $\cCa_{\leqa}$\,-\,operator and the $\cCb_{\leqa}$\,-\,operator.

\subsection{\textbf{$\cCa_{\leqa}$\,-\,operator}}\label{sec:agC1}

Recall that the $\cCa_{\leqa}$\,-\,operator has the form
\begin{align*}
    \cCa_{\leqa, \sigma, \mathrm{G}}(\mathbf{X},\mu)& = \bigoplus_{i=1}^n \rG(\bx_{\sigma(i)}, \mu(B_{\sigma(i)}),\mu(B_{\sigma(i+1)})),
\end{align*}
where $\rG\colon \sK \times [0,1]^2 \to \sK$ and $\sigma \in \Pi_{\bX}$.


\begin{corollary}\label{cor1:1}
Assume that $\lea$ is compatible with $\oplus.$
The $\cCa_{\leqa}$\,-\,operator is an~$n$-ary 
 $\sK$-valued aggregation function w.r.t.~$\leqa$ if and only if the following conditions hold:
\begin{itemize}[noitemsep]
    \item for  n=2: 
        \begin{enumerate}[label=(\alph*), noitemsep]
            \item 
             $[0,1] \ni c\mapsto \mathrm{G}(\bx, 1, c)\oplus\mathrm{G}(\bx, c, 0)$  is a~constant function  for any $\bx\in\sK$,
    
            \item  $\sK \ni \bx\mapsto \mathrm{G}(\bx, 1, b)$ and $\sK\ni \bx\mapsto \mathrm{G}(\bx, b, 0)$ are non-decreasing functions w.r.t.~$\leqa$ for any $b\in [0,1]$, 
    
            \item  $\mathrm{G}(\bzer, 1, b) \oplus\mathrm{G}(\bzer,b,0) = \bzer$ for any $b\in [0,1]$,
            
            \item $\mathrm{G}(\bjed, 1, b)\oplus\mathrm{G}(\bjed, b, 0) = \bjed$ for any $b\in [0,1]$.  
        \end{enumerate}

    \item for $n=3$: 
        \begin{enumerate}[label=(\alph*), noitemsep]
            \item 
            $[b, 1] \ni c \mapsto\rG(\bx, 1, c) \oplus \rG(\bx, c, b)$ and $[0, b] \ni c\mapsto \rG(\bx, b, c) \oplus \rG(\bx, c, 0)$ are constant functions for any $\bx \in \sK$ and $b\in [0,1]$,

            \item $\sK\ni \bx\mapsto \rG(\bx, b_1, b_2)$ is a~non-decreasing  function w.r.t.~$\leqa$ for any  $b_2\le b_1$, 
    
            \item  
            ${\tx \rG(\bzer, 1, b_1)\oplus  \rG(\bzer,b_1,b_{2})\oplus \rG(\bzer, b_2, 0)=\bzer}$ for any $b_2\le b_1$,

            \item ${\tx \rG(\bjed, 1, b_1)\oplus  \rG(\bjed, b_1, b_2)\oplus \rG(\bjed, b_2, 0)=\bjed}$ for any $b_2\le b_1$.

        \end{enumerate}

    \item for  $n\ge 4$:       
    \begin{enumerate}[label=(\alph*), noitemsep]
            \item 
            $[b_2, b_1] \ni c\mapsto \rG(\bx, b_1, c) \oplus \rG(\bx, c, b_2)$    is a~constant function  for any $\bx\in \sK$ and $b_2 \le b_1$, 
             
            \item $\sK \ni\bx\mapsto \rG(\bx, b_1, b_2)$ is a~non-decreasing  function w.r.t.~$\leqa$ for any $b_2\le b_1$, 
    
            \item ${\tx\bigoplus_{i=1}^n \rG(\bzer,b_i,b_{i+1}) = \bzer}$ for any $0=b_{n+1} \le b_n \le \ldots \le b_2 \le b_1=1$, 
            
            \item ${\tx \bigoplus_{i=1}^n \rG(\bjed,b_i,b_{i+1})=\bjed}$ for any $0=b_{n+1} \le b_n \le \ldots \le b_2 \le b_1=1$. 
        \end{enumerate}
\end{itemize}
\end{corollary}
\begin{pf}

Consider the case $n=2.$ Obviously, conditions (c) and (d) are equivalent to (b) and (c) of  Theorem~\ref{tw:ag}.
In view of Remark~\ref{rem:3.4}, condition (a) is equivalent to condition \ref{WD2} of Theorem \ref{tw:well}. 
From Definition~\ref{def:2.4} and Lemma~\ref{lem:3}\,(b) 
it follows that condition (b) corresponds to conditions (b) and (c) of Theorem~\ref{tw:mon_2}.
The statements for $n = 3$ and $n\ge 4$ can be proved similarly.
\qed
\end{pf}

Putting $\sK = [0,1]$ in Corollary~\ref{cor1:1}, we get a~refinement for $n=3$ of Theorem~2.4 from \cite{BJKO23}.
In the next two sections, 
namely \ref{sec:4.1.1} and \ref{sec:Gdis}, we will look at two special forms of the function $\rG$.

\subsubsection{$\rG(\bx, b_1, b_2) = \delta(b_1,b_2) \odot \bx$}
\label{sec:Gdis}

 In this section, we assume that $\rG(\bx, b_1, b_2) =  \delta(b_1,b_2) \odot \bx$ for any $\bx \in \sK$ and $b_2 \le b_1$,  where $\delta$ is a~real-valued dissimilarity function (see Definition~\ref{def:dis}) and
 $\odot \colon [0,1] \times \sK \to \sU$
is an operation, called the \textit{multiplication},  satisfying the following conditions:
\begin{enumerate}[noitemsep, label=(\alph*)]
    \item 
    $1\odot \bx = \bx$ for any $\bx$, 
    
    \item $c\odot \bzer = \bzer$ for any $c$. 
\end{enumerate}   

 Examples of multiplication operations for the universal sets discussed in  Example~\ref{ex:1} are  given in Table~\ref{tab:2}.

\begin{table}[h!]
\centering
\begin{tabular}{|c|c| c| c|c|}\hline
 $\sU$  & $\sK$ &  $\odot$ & $c \odot \bz$ \\ \hline
$\mR_+$  & $[0,1]$ & $\cdot$  & $c\cdot z=cz$ (standard multiplication) \\ \hline
 $\int$   & $\sLi$ & $\odot_{\mathrm{iv}}$ & $c\odot_{\mathrm{iv}} [z^l,z^u]=[cz^l, cz^u]$ \\ \hline
$\mR^k_+$  & $[0,1]^k$ &$\odot_{\mathrm{vv}}$ & $c\odot_{\mathrm{vv}}(z_1,\ldots,z_k)^T=(cz_1,\ldots, cz_k)^T$  \\ \hline
\end{tabular}
    \caption{Examples of  multiplication  operations  for  selected  cases of $\sU$.}
    \label{tab:2}
\end{table}

Before we present a~characterization of all real-valued dissimilarity functions for which the $\cCa_{\leqa}$\,-\,operator is a~$\sK$-valued aggregation function, let us introduce some  auxiliary concepts.

\begin{definition}\label{def:4.9}
    Let $\odot$ be a~multiplication operation. 
    \begin{itemize}[noitemsep]
        \item A partial order $\leqp$  is said to be \textit{compatible   with $\odot$} if for all $c\in [0,1]$, $c\odot \bx \leqp c\odot \by$ whenever $\bx\leqp \by$.

        \item The operation $\odot$ has the \textit{cancellation law}, whenever for any $c_1,c_2 \in [0,1]$ the condition $c_1 \odot \bjed = c_2 \odot \bjed$ implies $c_1 = c_2$.

        \item  The operation $\odot$ is \textit{right-distributive} over $\oplus$ if $(c_1 + c_2) \odot \bx = (c_1 \odot \bx) \oplus (c_2 \odot \bx)$ for any $\bx\in \sK$ and  $c_1,c_2\in [0,1]$.
    \end{itemize}
\end{definition}

\begin{corollary}\label{cor:4.5}
Let $n\ge 3$, 
and 
a~multiplication operation $\odot$ is right-distributive over $\oplus$ and has the cancellation law.
Assume that
$\lea$ is compatible with $\oplus$ and $\leqa$ is compatible with $\odot$.
The $\cCa_{\leqa}$\,-\,operator is an~$n$-ary $\sK$-valued aggregation function w.r.t.~$\leqa$ if and only if 
\begin{align}\label{rdes:1}
\delta(b_1, b_2) = \delta(b_1,0) - \delta( b_2, 0) \quad \text{ for any } b_2\le b_1.
\end{align}
\end{corollary}
\begin{pf}
Let us begin with the case $n\ge 4$.
The condition (b) of Corollary~\ref{cor1:1} follows from the fact that $\leqa$ is compatible with $\odot$.


Due to the right-distributivity of $\odot$, condition (a)  of Corollary~\ref{cor1:1} takes the form
\begin{align}\label{pomdis:0}
    (\delta(b_1, c) + \delta(c,b_2)) \odot \bx = (\delta(b_1, d) + \delta(d, b_2))\odot \bx
\end{align}
for any $\bx\in \sK$, $b_2 \le c\le b_1$, and $b_2 \le d\le  b_1$.
Putting  $\bx = \bjed$ in \eqref{pomdis:0}  and using the cancellation law  of $\odot$, we can rewrite condition (a) of Corollary~\ref{cor1:1} as 
\begin{align}\label{pomdis:1}
\delta(b_1, c) + \delta(c,b_2) = \delta(b_1, d) + \delta(d, b_2)\quad \text{ for any }  b_2 \le c\le b_1 \text{ and } b_2 \le d\le  b_1.
\end{align}
Taking $b_2=0=c$ yields
\begin{align}\label{pomdis:2}
    \delta(b,d) = \delta(b,0) - \delta(d,0) \quad \text{ for any } d\le b.
\end{align}
It is clear that \eqref{pomdis:2} implies \eqref{pomdis:1}, so these conditions are equivalent to each other.


Because of the right-distributivity of $\odot$, condition (c)  of Corollary~\ref{cor1:1} has form  
\begin{align*}
    \big(\sum_{i=1}^n \delta(b_i, b_{i+1})\big)\odot \bzer=\bzer \quad \text{ for any } 0=b_{n+1} \le b_n\le \ldots \le b_2 \le b_1 =1.
\end{align*}
By \eqref{pomdis:2} and Definition~\ref{def:dis},
${\tx \sum_{i=1}^n \delta(b_i, b_{i+1})=\delta(1,0)=1}$. 
From the definition of $\odot$, it follows that $1\odot \bzer = \bzer$, so 
condition (c)  of Corollary~\ref{cor1:1} is valid.
Similarly, it can be shown that condition (d) of Corollary~\ref{cor1:1} also holds. 
The proof for $n=3$ is similar, so we will omit it. 
\qed
\end{pf}

As a~direct consequence of Corollary~\ref{cor:4.5} we get the following result for $\sK = \sLi$.

\begin{corollary}
Let $n\ge 3$ and $\sK = \sLi$. 
Assume that $\lea$ is compatible with $\oplus_{\mathrm{iv}}$ and $\leqa$ is compatible with  $\odot_{\mathrm{iv}}$.
The $\cCa_{\leqa}$\,-\,operator
$$
\cCa_{\leqa, \sigma,\delta}(\bX, \mu) = \bigoplus_{i=1}^{n} \!{}_{\tiny\mathrm{iv}} \big(\delta(\mu(B_{\sigma(i)}), \mu(B_{\sigma(i+1)})) \odot_{\mathrm{iv}} \bx_{\sigma(i)}\big),
$$ 
where
$\bX,\mu,\sigma, B_{\sigma(i)}$ are as in~\eqref{dc}, is an~$n$-ary IV aggregation function w.r.t.~$\leqa$ if and only if  \eqref{rdes:1}  holds. 
\end{corollary}

\begin{remark}\label{rem:4.7a}
Clearly, $\lea_{(\alpha,\beta)}$ is compatible with $\oplus_{\mathrm{iv}}$ and $\leqab$ is compatible with $\odot_{\mathrm{iv}}$. 
Since $\oplus_{\mathrm{iv}}$ has the cancellation law, due to Lemmas~\ref{lem:1}\,-\,\ref{lem:2}, $\leqa$ is compatible with $\oplus_{\mathrm{iv}}$, and $\lea$ is compatible with $\oplus_{\mathrm{iv}}$ whenever at least one of them is. 
By \cite[Cor.~2.7]{BJK21}, the only admissible order that is compatible with $\oplus_{\mathrm{iv}}$ and $\odot_{\mathrm{iv}}$ is the $(\alpha,\beta)$-order.  
Hence, if  $\delta(b_1,b_2) = b_1 - b_2$ for $b_2\le b_1$, then the only functions 
$A,B\colon [0,1]^2\to [0,1]$ 
for which the $\preceq_{A,B}$\,-\,Choquet integral\footnote{see \cite[Def.~4]{bustince2013b}}  is an IV aggregation function w.r.t.~$\preceq_{A,B}$ within the class of 
admissible orders  $\preceq_{A,B}$ compatible with
$\oplus_{\mathrm{iv}}$ and $\odot_{\mathrm{iv}}$  are $A(x, y)= (1-\alpha) x + \alpha y$ and $B(x, y) =(1-\beta) x + \beta y $ with $\alpha,\beta\in [0,1]$ and $\alpha\neq \beta$.
\end{remark}

We conclude this section with the following direct consequence of Corollary~\ref{cor:4.5} for $\sK=[0,1]^k$.

\begin{corollary}
Let $n\ge 3$ and $\sK = [0,1]^k$ with $k\in \mN$.
Assume that $\lea$ is compatible with $\oplus_{\mathrm{vv}}$ and $\leqa$ is compatible with $\odot = \odot_{\mathrm{vv}}$.
The $\cCa_{\leqa}$\,-\,operator is an $n$-ary  VV aggregation function w.r.t.~$\leqa$ if and only if \eqref{rdes:1} holds.
\end{corollary}

For $k=1$ we get Corollary~2.7 from \cite{BJKO23}.


\subsubsection{$\rG(\bx, b_1, b_2) = \rF(\bx, b_1 - b_2)$}\label{sec:4.1.1}

In this part, we assume that $\rG(\bx, b_1, b_2) = \rF(\bx, b_1 - b_2)$ for any $\bx \in \sK$ and $b_1 \ge b_2$ with $\rF\colon \sK \times [0,1] \to \sK$. 
From  Corollary~\ref{cor1:1} for $n\ge 3$ and Remark~\ref{rem:42} we obtain the following result.

\begin{corollary}\label{cor:4.9}
Assume that $n\ge 3$,  $\bzer \oplus \bzer \in \sK$, and $\lea$ is compatible with $\oplus$.
The $\cCa_{\leqa}$\,-\,operator is an~$n$-ary $\sK$-valued aggregation function 
w.r.t.~$\leqa$ if and only if the following conditions hold:
\begin{enumerate}[noitemsep, label=(\alph*)]
    \item $\rF(\bx, a) \oplus \rF(\bx, b) = \rF(\bx, 0.5(a+b)) \oplus \rF(\bx, 0.5(a+b))$  for all $\bx \in \sK$ and $a,b \in [0,1]$,

    \item $\sK \ni\bx\mapsto \rF(\bx, b)$ is a~non-decreasing  function w.r.t.~$\leqa$ for any  $b\in [0,1]$,

    \item $\rF(\bzer, b) = \bzer$ for any $b\in [0,1]$,
    
    \item $\tx\bigoplus_{i=1}^n \rF(\bjed, b_i) = \bjed$ for all $b_1,\ldots,b_n \in [0,1]$ such that $b_1+\ldots+b_n=1$. 
\end{enumerate}  
\end{corollary}

Let $\sK =\sLi$.
Then there exist $\rF^l,\rF^u\colon \sLi \times [0,1] \to [0,1]$  such that $\rF^l\le \rF^u$ and
 $\rF(\bx, b)= [\rF^l(\bx, b), \rF^u(\bx, b)].$
Thus Corollary~\ref{cor:4.9}\,(a) takes the form 
of two Jensen equations:
$$
\rF^i(\bx,a)+\rF^i(\bx,b)=2\rF^i(\bx, 0.5 (a+b)), \quad i\in \{l,u\},
$$
for all $\bx\in \sLi$ and $a,b\in [0,1]$. 
Since 
$c\mapsto \rF^i(\bx, c)$
is bounded for each $\bx$, by \cite[Sec.~1]{kanapan2009}, we have
$$
\rF^i(\bx,a) = aC^i(\bx) + D^i(\bx), \qquad i\in \{l,u\},
$$
for all $\bx\in\sLi$, $a\in [0,1]$, and some $C^i,D^i \colon \sLi \to [0,1]$ such that $a C^l(\bx)+D^l(\bx) \le aC^u(\bx)+D^u(\bx)$.
Combining  the above with points (b)-(d)  of  Corollary~\ref{cor:4.9}, we get the following  result.

\begin{corollary}
Let $n\ge 3$ and $\sK = \sLi$. 
Assume that $\lea$ is compatible with $\oplus_{\mathrm{iv}}$.
The $\cCa_{\leqa}$\,-\,operator is an $n$-ary IV aggregation function w.r.t.~$\leqa$ if and only if
\begin{align*}
\rF(\bx,a) = [aC^l(\bx) + D^l(\bx),\, aC^u(\bx) + D^u(\bx)]
\end{align*}
for all $\bx \in \sLi$ and $a\in [0,1]$ with some functions $C^i,D^i\colon \sLi \to [0,1]$ for $i\in\{l,u\}$ such that:
\begin{itemize}[noitemsep]
    \item $a C^l(\bx) + D^l(\bx) \le a C^u(\bx) + D^u(\bx)$  for any $\bx\in \sLi$  and $a\in [0,1]$,
    
   \item   $\sLi \ni\bx \mapsto aC^i(\bx) + D^i(\bx)$, $i\in \{l,u\},$ are non-decreasing functions w.r.t.~$\leqa$,
    
    \item $C^i(\bzer)=D^i(\bzer)=0$  and $C^i(\bjed) + n D^i(\bjed)=1$ for $i\in \{l,u\}$, where $\bzer = [0,0]$ and $\bjed = [1,1]$.
\end{itemize}
\end{corollary}

\medskip

Using the same arguments  as for $\sK = \sLi$, we  obtain from Corollary~\ref{cor:4.9} the result  for $\sK = [0,1]^k$.

\begin{corollary}
Let $n\ge 3$ and $\sK = [0,1]^k$ with $k\in \mN$.
Assume that $\lea$ is compatible with $\oplus_{\mathrm{vv}}$.
The $\cCa_{\leqa}$\,-\,operator is an $n$-ary  VV aggregation function w.r.t.~$\leqa$ if and only if
\begin{align*}
\rF(\bx,a) = (aC_1(\bx) + D_1(\bx), \ldots, aC_k(\bx) + D_k(\bx))^T
\end{align*}
for all $\bx\in [0,1]^k$ and $a\in [0,1]$ with some functions $C_i,D_i\colon [0,1]^k \to [0,1]$ such that:
\begin{itemize}[noitemsep]
    \item   $[0,1]^k \ni\bx \mapsto aC_i(\bx) + D_i(\bx)$, $i\in [k]$,  are non-decreasing functions w.r.t.~$\leqa$,
    
    \item $C_i(\bzer)=D_i(\bzer)=0$ and $C_i(\bjed) + n D_i(\bjed)=1$ for $i\in [k]$ with $\bzer = (0,\ldots,0)^T$ and $\bjed = (1,\ldots,1)^T$.
\end{itemize}
\end{corollary}

 Note that for $k=1$ we recover Corollary~2.6 from \cite{BJKO23}.

\subsection{\textbf{$\cCb_{\leqa}$\,-\,operator}}\label{sec:agC2}

We will now examine the conditions that characterize the $\cCb_{\leqa}$\,-\,operator as a~$\sK$-valued aggregation function.
Recall that the $\cCb_{\leqa}$\,-\,operator  is of the form
\begin{align*}
    \cCb_{\leqa, \sigma, \mathrm{G}}(\mathbf{X},\mu)& = \bigoplus_{i=1}^n \mathrm{G}(\bx_{\sigma(i)}, \bx_{\sigma(i-1)}, \mu(B_{\sigma(i)})),
\end{align*}
where $\rG\colon \sK^2 \times [0,1] \to \sK$ and $\sigma \in \Pi_{\bX}$.

\begin{corollary}\label{cor:4.12}
Assume that  $\lea$ is compatible with $\oplus.$
The $\cCb_{\leqa}$\,-\,operator is an~$n$-ary  $\sK$-valued aggregation function w.r.t.~$\leqa$ if and only if the following conditions are met:
\begin{itemize}
    \item for  $n=2$: 
        \begin{enumerate}[label=(\alph*), noitemsep]
             \item $\rG(\bx, \bx, b) =  \rG(\bx, \bx, 0)$ for any $\bx\in \sK$ and $b\in [0,1]$, 

            \item  $[\bzer, \bv]_{\sK}\ni\bx \mapsto\rG(\bx, \bzer, 1) \oplus \mathrm{G}(\bv, \bx, b)$ 
            is a~non-decreasing function w.r.t.~$\leqa$ for any $\bv\in \sK$ and $b \in [0,1]$,

            \item 
            $[\bu, \bjed]_{\sK} \ni \bx \mapsto \mathrm{G}(\bx, \bu, b)$
            is a~non-decreasing function w.r.t.~$\leqa$ for any $\bu\in \sK$ and $b \in [0,1]$,

            \item $\rG(\bzer,\bzer,0) \oplus \rG(\bzer,\bzer,0) = \bzer$,
            
            \item  $\rG(\bjed,\bzer,1)\oplus\rG(\bjed,\bjed,0)=\bjed$. 
        \end{enumerate}

    \item for  $n\ge 3$: 
        \begin{enumerate}[label=(\alph*), noitemsep]
            \item   $\rG(\bx, \bx, b) = \rG(\bx, \bx, 0)$ for any $\bx\in \sK$ and $b\in [0,1]$, 
           
            \item $[\bu, \bv]_{\sK}\ni \bx\mapsto 
            \rG(\bx, \bu, b_1) \oplus \rG(\bv, \bx,  b_2)$ 
            is a~non-decreasing  function w.r.t.~$\leqa$ for any 
            $\bu\leqa \bv$ and $b_2 \le b_1$,
    
            \item $[\bu, \bjed]_{\sK}\ni\bx\mapsto \rG(\bx, \bu, b)$ 
            is a~non-decreasing function w.r.t.~$\leqa$ for any $\bu\in \sK$ and $b \in [0,1]$,

            \item  ${\tx\bigoplus_{i=1}^n \rG(\bzer,\bzer, 0) = \bzer}$, 
            
            \item ${\tx \rG(\bjed,\bzer,1)\oplus \bigoplus_{i=2}^n \rG(\bjed,\bjed,0)=\bjed}$.
        \end{enumerate}
\end{itemize}
\end{corollary}
\begin{pf}
To prove these statements, we use Theorem~\ref{tw:ag} with $\rL(\bx_1,\bx_2, b_1, b_2)=\rG(\bx_1, \bx_2, b_1)$.
Let $n=2.$
It follows from 
Remark~\ref{rem:3.4}
that condition (a) is equivalent to 
condition (a) of Theorem~\ref{tw:mon_2}.
Conditions (b) and (c) represent conditions (b) and (c) of Theorem~\ref{tw:mon_2}, while, by (a), conditions (d) and (e) are equivalent to conditions (b) and (c) of Theorem~\ref{tw:ag}. 
Similarly, one can prove the statements for $n=3$ and $n\ge 4$, and observe that they are identical. \qed
\end{pf}

\subsubsection{$\rG(\bx_1, \bx_2, b) = b\odot \rd(\bx_1, \bx_2)$}\label{Sec:4.2.1}

We now consider $\rG(\bx_1, \bx_2, b) = b\odot \rd(\bx_1, \bx_2)$ for $\bx_2 \leqa \bx_1$,  where $\rd$ is a~$\sK$-valued dissimilarity function  (see Definition~\ref{def:dis}) and $\odot$ is a~multiplication operation  (see Section~\ref{sec:Gdis}). 
 We will need the following assumption
\begin{enumerate}
    \item[(C1)] $(b_1 \odot \bu_1) \oplus (b_2\odot \bv_2) \leqa (b_1 \odot\bu_2) \oplus (b_2 \odot\bv_1)$ for any $b_i\in [0,1]$ and $\bu_i,\bv_i\in \sK$ such that $b_2\le b_1$,
    $\bu_1\leqa \bu_2$, $\bv_1\leqa \bv_2$, and $\bu_1 \oplus \bv_2 = \bu_2 \oplus \bv_1\in\sK$.
\end{enumerate}

We say that a~multiplication operation $\odot$ is \textit{left-distributive over $\oplus$} if $c \odot (\bx \oplus \bz) = (c\odot \bx ) \oplus (c\odot \bz)$ for any $\bx,\bz\in \sK$ and $c\in [0,1]$ such that $\bx \oplus \bz \in \sK$.

\begin{corollary}\label{cor:4.13}
Let $n\ge 3$, 
$\odot$ be left-distributive over $\oplus$, and $\bzer\oplus \bx = \bx$ for any  $\bx\in \sK$. 
Assume that (C1) holds,  $\lea$ is compatible with $\oplus$ and $\leqa$ is compatible with $\odot$.
The $\cCb_{\leqa}$\,-\,operator is an~$n$-ary $\sK$-valued aggregation function w.r.t.~$\leqa$ if and only if
\begin{align}\label{dis:2}
    \rd(\bx_1, \bzer) \oplus\rd(\bx_2, \bx_1) = \rd(\bx_2, \bzer)\qquad \text{ for any }  \bx_1,\bx_2\in \sK \text{ such that } \bx_1 \leqa \bx_2.
\end{align}
\end{corollary}
\begin{pf}
Conditions (a), (d), and (e)  of  Corollary~\ref{cor:4.12} are satisfied in view of definitions of the $\sK$-valued dissimilarity function and the operation $\odot$. 
Point (c) of  Corollary~\ref{cor:4.12} follows from Definition~\ref{def:dis}\,(d) and the fact that $\leqa$ is compatible with $\odot$.
Point (b) of Corollary~\ref{cor:4.12} takes the form
\begin{align}\label{dis:1}
 (b_1 \odot \rd(\bx_1, \bu)) \oplus (b_2 \odot\rd(\bv, \bx_1)) \leqa  (b_1 \odot\rd(\bx_2, \bu)) \oplus (b_2 \odot\rd(\bv, \bx_2))
\end{align}
for any $b_2\le b_1$ and 
$\bx_1 \leqa \bx_2$ such that $\bx_1,\bx_2 \in [\bu,\bv]_{\sK}$.
We now show that \eqref{dis:1} 
is equivalent to \eqref{dis:2}.

``$\Rightarrow$''
Putting $b_1=b_2=1$ in \eqref{dis:1} and using the fact that $1 \odot \bx = \bx$ for any $\bx$ we have
\begin{align}\label{dis:1a}
    \rd(\bx_1, \bu) \oplus \rd(\bv, \bx_1) \leqa  \rd(\bx_2, \bu) \oplus \rd(\bv, \bx_2)
\end{align}
for any $\bx_1 \leqa \bx_2$ such that $\bx_1,\bx_2 \in [\bu,\bv]_{\sK}$.
Putting $\bu = \bx_1 = \bzer$ in \eqref{dis:1a} gives $\rd(\bv, \bzer) \leqa \rd (\bx_2, \bzer)\oplus \rd (\bv, \bx_2)$ 
for any 
$\bx_2 \in [\bzer, \bv]_{\sK}$,
as $\bzer \oplus \bx = \bx$ for any $\bx\in \sK$.
Setting $\bx_2 = \bv$ and $\bu = \bzer$ in \eqref{dis:1a} yields $\rd(\bx_1, \bzer) \oplus \rd(\bv, \bx_1) \leqa \rd(\bv, \bzer)$
for any $\bx_1 \in [\bzer, \bv]_{\sK}$,
as  $\bx \oplus \bzer = \bx$ for any $\bx\in \sK$.
So by the transitivity of $\leqa$, we get \eqref{dis:2}.

``$\Leftarrow$''
From \eqref{dis:2}  we have
\begin{align}
\rd (\bx_1, \bzer) \oplus \rd(\bv, \bx_1) = \rd (\bx_2, \bzer) \oplus \rd(\bv, \bx_2)\label{mb2}
\end{align}
for any 
$\bx_1 \leqa \bx_2$ such that $\bx_1,\bx_2 \in [\bzer,\bv]_{\sK}$.
 Additionally, from \eqref{dis:2} we know that the left and right hand side of \eqref{mb2} is an element of $\sK$.
Since $\rd (\bx_1, \bzer)\leqa \rd (\bx_2, \bzer)$ and $\rd(\bv, \bx_2)\leqa \rd(\bv, \bx_1)$ for any $\bx_1\leqa \bx_2 \leqa \bv$, so by (C1) and \eqref{mb2}, we have
 $L\leqa P$,
where $L= (b_1 \odot \rd(\bx_1, \bzer))  \oplus (b_2\odot \rd(\bv, \bx_1))$
and $P= (b_1 \odot \rd(\bx_2, \bzer))  \oplus (b_2\odot \rd(\bv, \bx_2))$
for any $b_2\le b_1$ and  $\bx_1 \leqa \bx_2$ such that $\bx_1,\bx_2 \in [\bzer, \bv]_{\sK}$.
 Applying \eqref{dis:2}, we  obtain
\begin{align}\label{mb3}
    \rd(\bx_1, \bzer)= \rd(\bu, \bzer) \oplus \rd(\bx_1, \bu), \qquad
    \rd(\bx_2, \bzer)= \rd(\bu, \bzer) \oplus \rd(\bx_2, \bu)
\end{align}
for any $\bx_1 \leqa \bx_2$ such that $\bx_1,\bx_2 \in [\bu, \bjed]_{\sK}$.
Using \eqref{mb3} and the left-distributivity of $\odot$,  
we get 
\begin{align}
L&=(b_1\odot \rd(\bu, \bzer)) \oplus (b_1 \odot\rd(\bx_1, \bu)) \oplus (b_2 \odot\rd(\bv, \bx_1))\label{dis:4a},\\
P&=  (b_1 \odot\rd(\bu, \bzer)) \oplus (b_1 \odot\rd(\bx_2, \bu)) \oplus (b_2\odot\rd(\bv, \bx_2))  \label{dis:4b}
\end{align}
for any  $b_2\le b_1$ and 
$\bx_1 \leqa \bx_2$ such that $\bx_1,\bx_2 \in [\bu,\bv]_{\sK}$. 
Combining \eqref{dis:4a}\,-\,\eqref{dis:4b}
with Lemma~\ref{lem:3}\,(b)
we obtain \eqref{dis:1} for any $ b_2\le b_1$ and
$\bx_1 \leqa \bx_2$ such that $\bx_1,\bx_2 \in [\bu,\bv]_{\sK}$. 
\qed
\end{pf}

If $\sK = [0,1]$, then we  recover Corollary~3.9 from \cite{BJKO23}. The case when $\sK = \sLi$  and the admissible order is the $(\alpha,\beta)$-order (see Remark~\ref{rem:4.7a})
is summarized below.

\begin{corollary}\label{cor:4.14}
Let $n\ge 3$ and $\sK = \sLi$.
The following conditions are equivalent.
\begin{enumerate}[noitemsep, label = (\alph*)]
    \item The $\cCb_{\leqab}$\,-\,operator is an $n$-ary IV aggregation function w.r.t.~$\leqab$.
    
    \item  $\rd(\bx_1, \bzer) \oplus_{\mathrm{iv}} \rd(\bx_2, \bx_1) = \rd(\bx_2, \bzer)$ for any $\bx_1 \preceq_{(\alpha, \beta)} \bx_2$.

    \item For each $k\in [n]$ it holds
    \begin{align*}
        \bigoplus_{i=1}^k\rd(\bx_{i-1}, \bx_i) \leqab \bigoplus_{i=1}^k \rd(\bz_{i-1}, \bz_i)
    \end{align*}
    for all $\bX = (\bx_1, \ldots, \bx_n)\in\sK^n$ and $\bZ=(\bz_1,\ldots, \bz_n)\in\sK^n$ such that $\bx_1 \leqab\ldots \leqab \bx_k$, $\bz_1 \leqab\ldots \leqab \bz_k$ and $\bx_1 \leqab \bz_1,\ldots, \bx_k \leqab \bz_k$ with $\bx_0=\bz_0=\bzer$, where here ${\tx \bigoplus_{i=1}^k \bx_i} = \bx_1 \oplus_{\mathrm{iv}} \ldots\oplus_{\mathrm{iv}} \bx_k$. 
\end{enumerate}
\end{corollary}
\begin{pf}
The equivalence of conditions (a) and (b) follows from Corollary~\ref{cor:4.13} with $\sK=\sLi$.
The equivalence of (a) and (c)  can be established  using Proposition~5.3 and Theorem~5.5 of  \cite{takac2022}. \qed
\end{pf}

The authors of \cite{takac2022} claim that it is a difficult task to find an IV dissimilarity function $\rd$ that satisfies point (c) of Corollary~\ref{cor:4.14}. 
In Appendix C, we show that there is no such function in a~certain subclass of IV dissimilarity functions considered in \cite[Prop.~4.18]{takac2022}.


\section{Conclusion}
In this article, we have introduced the Choquet-like operator, which is a~generalization of the classical discrete Choquet integral to a~multivalued setting in terms of a~general function that replaces subtraction and multiplication, and an admissible order that refines the natural partial order on the data set under consideration. The new operator takes an input consisting of a~finite number of values from the data set and returns a~single output value belonging to that set. We have characterized the class of all functions for which the Choquet-like operator is invariant with respect to permutations of the input data, and then provided necessary and sufficient conditions for the operator to exhibit monotonicity with respect to the admissible order, as well as a~complete characterization of the operator as an aggregation function with respect to the admissible order. 
We have also studied two important special types of the Choquet-like operator and a~number of their particular forms, both considered and not considered in the literature, for data of general type as well as for scalar, interval, and vector data. 
In particular, we have shown that the only order 
in a~certain class of admissible orders for which the $\preceq_{A,B}$\,-\,Choquet integral defined in \cite[Def.~4]{bustince2013b} is an interval-valued aggregation function with respect to $\preceq_{A,B}$, is the $(\alpha,\beta)$-order.

\section*{Declaration of interests}
The authors declare that they have no known competing financial interests or
personal relationships that could have appeared to influence the work reported
in this paper.

\section*{Data availability}

No data was used for the research described in the article.



\section*{Appendix A}

To shorten notation, from now on we will write $\cC_{\sigma}(\mathbf{x},\mu)$ instead of $\cC_{\leqa, \sigma, \rL}(\mathbf{x},\mu)$.

\begin{pot1}
``$\Rightarrow$''
Let $\mu\in \bM$.
Take $\bX= (\bx_1,\ldots,\bx_n)\in \sK^n$
 such that   $\bx_m=\bx_l=\bx$ for some distinct $m,l\in [n]$. 
Then there are two distinct admissible permutations w.r.t.~$\leqa$ for $\bX$, say $\sigma_1$ and $\sigma_2$, which satisfy for some $k\in \{2,\ldots,n\}$ the condition that $\bx_{\sigma_1(k)}=\bx_m$, $\bx_{\sigma_1(k-1)}=\bx_l$, $\bx_{\sigma_2(k-1)}=\bx_m$, $\bx_{\sigma_2(k)}=\bx_l$ and ${\sigma_1(i)}={\sigma_2(i)}$ for any $i\in [n]\setminus \{k-1,k\}$. 
Clearly, 
 $B_{\sigma _1(i)}=B_{\sigma _2(i)}$ 
for each $i\in [n]\setminus \{k\}$. 
Since, by the assumption,  $\cC_{\sigma _1}(\bX,\mu) = \cC_{\sigma_2}(\bX,\mu)$, applying  Definition~\ref{def:choquet} and 
the cancellation law of $\oplus$, we get 
\begin{align*}
    &\rL(\bx, \bx_{\sigma_1(k-2)}, \mu(B_{\sigma _1(k-1)}),\mu(B_{\sigma _1(k)})) \oplus  \rL(\bx, \bx, \mu(B_{\sigma _1(k)}),\mu(B_{\sigma _1(k+1)}))\notag
    \\&\quad =
   \rL(\bx, \bx_{\sigma_1(k-2)}, \mu(B_{\sigma _1(k-1)}),\mu(B_{\sigma_2(k)})) \oplus  \rL(\bx, \bx, \mu(B_{\sigma_2(k)}),\mu(B_{\sigma _1(k+1)})).
\end{align*}
Set $b_i=\mu(B_{\sigma _1(i)})$ and $d_i=\mu(B_{\sigma_2(i)})$ for any $i\in [n+1]$  and $\widehat{\bx}_{k-2}=\bx_{\sigma_1(k-2)}$. 
Clearly, $b_i = d_i$ for any $i\in [n+1]\setminus\{k\}$, 
$b_{k+1}\le b_k \le b_{k-1}$, and $b_{k+1}\le d_k \le b_{k-1}$.
 Hence, for $k\in \{2,\ldots ,n\}$,
\begin{align}\label{cc4}
\rL(\bx, \widehat{\bx}_{k-2}, b_{k-1}, b_k) \oplus  \rL(\bx, \bx, b_k,b_{k+1})  =
\rL(\bx,\widehat{\bx}_{k-2}, b_{k-1}, d_k) \oplus \rL(\bx, \bx, d_k, b_{k+1}), 
\end{align}
and consequently
 $    [b_{k+1}, b_{k-1}]\ni c \mapsto \rL(\bx, \widehat{\bx}_{k-2}, b_{k-1}, c) \oplus  \rL(\bx, \bx, c, b_{k+1})$ is a~constant function.
Due to the arbitrariness of $\mu$, this leads us to 
\ref{WD2} for $n=2$, \ref{WD3} for $n=3$ and \ref{WDn} for $n\ge 4$.  

``$\Leftarrow$'' 
Let $\bX= (\bx_1,\ldots,\bx_n)\in \sK^n$ be such that there are 
two  permutations $\sigma_1, \sigma_2 \in \Pi_{\bX}$,  $\sigma_1 \not= \sigma_2. $  
This means that there are at least two elements 
in $\bX$ that have the same value, and 
$\sigma_2$ is a~permutation that swaps equal elements in the sequence $(\bx_{\sigma_1(1)}, \ldots, \bx_{\sigma_{1}(n)})$.  
It is evident that  
$\sigma _2$ can be represented as the product of $\sigma_1$ and some 
adjacent transpositions (see \cite[p.~631]{olver2021}), say, $\tau_1,\ldots, \tau_m$ with certain $m\ge 1$, that is, 
$\sigma _2 =\tau_m\circ \ldots \circ \tau_1 \circ \sigma _1$.  
So, for any $j\in [m]$, 
$\tau_{j}$  is the adjacent  transposition that exchanges the elements $\bx_{\pi_{j-1}(k_j-1)}$ and  $\bx_{\pi_{j-1}(k_j)}$ equal to each other for some $k_j\in \{2,3,\ldots ,n\}$ in the sequence 
$(\bx_{\pi_{j-1}(1)}, \ldots, \bx_{\pi_{j-1}(n)})$, where 
$\pi_j=\tau_j\circ \pi_{j-1}$ with $\pi_0 = \sigma_1$.
Clearly, 
$\pi_j \in\Pi_{\bX}$ for $j\in [m]$.
By \ref{WD2}, \ref{WD3} or \ref{WDn}, we get \eqref{cc4} with  $\bx=\bx_{\pi_j(k_j-1)}$, $\widehat{\bx}_{k-2}=\bx_{\pi_j(k_j-2)}$,  $b_{k-1}=\mu(B_{\pi_j(k_j-1)})$, 
$b_k=\mu(B_{\pi_j(k_j)})$, $b_{k+1}=\mu(B_{\pi_j(k_j+1)})$, and  $d_k=\mu(B_{\pi_{j+1}(k_j)})$ 
for $j\in \{0,1,\ldots, m-1\}$,
which implies that $\cC_{\pi_j}(\bX,\mu) = \cC_{\pi_{j+1}}(\bX,\mu)$ for $j\in \{0,1,\ldots, m-1\}$. 
Hence $\cC_{\sigma_1}(\bX,\mu) = \cC_{\sigma _2}(\bX,\mu)$. 
This completes the proof.
\qed
\end{pot1}

\section*{Appendix B}

Below, we will present a~proof of the monotonicity characterization for the Choquet-like operator. We will use the convention outlined in Appendix A, and denote by $\mathrm{id}$ the identity function on $[n]$. We shall need an auxiliary result using the following condition:

\begin{enumerate}[label=(M'), leftmargin=1.4cm]
    \item \label{ma} 
    Assume that $\bX=(\bx_1,\ldots,\bx_n)\in \sK^n$ and $\bZ=(\bz_1,\ldots,\bz_n)\in\sK^n$ are such that $\bx_1\leqa \ldots\leqa \bx_n$,  $\bx_p\lea \bz_p$  for some $p\in [n]$ and $\bx_i = \bz_i$ for any $i\in [n]\setminus \{p\}$.
    Then
    $\cC_{\mathrm{id}}(\bX,\mu)\leqa \cC_{\tau}(\bZ,\mu)$ for any $\mu\in \bM$ and some $\tau \in \Pi_\bZ$. 
\end{enumerate}

\begin{lemma}\label{lem:10}
Assume that the $\cC_{\leqa}$\,-\,operator is well defined.
If condition~\ref{ma} is true, then condition~\ref{m1} holds.  
\end{lemma}
\begin{pf}
Let $\bX = (\bx_1,\ldots,\bx_n)\in\sK^n$ and 
$\bZ=(\bz_1,\ldots,\bz_n)\in\sK^n$ be such that $\bx_p \lea \bz_p$ for some $p\in [n]$ and $\bx_i=\bz_i$ for any $i\in [n]\setminus \{p\}$, so $\bX\leqA \bZ$. 
Firstly, we show  that  
\begin{align}\label{n:20}
    \cC_{\sigma}(\bX,\mu)\leqa \cC_{\tau}(\bZ,\mu)\qquad \text{ for any } \mu\in \bM, \text{ } \sigma\in \Pi_{\bX}, \text{ and } \tau\in \Pi_{\bZ}.
\end{align}
Fix $\mu\in \bM$, $\sigma\in \Pi_{\bX}$, and $\tau\in \Pi_{\bZ}$.
Define $\widehat{\mu}\in \bM$ as 
$\widehat{\mu}(C)=\mu(\{\sigma(i)\mid i\in C\})$ for any $C\subseteq [n]$. 
Let $\widehat{\bX}=(\widehat{\bx}_1,\ldots,\widehat{\bx}_n)\in \sK^n$ and $\widehat{\bZ}=(\widehat{\bz}_1,\ldots,\widehat{\bz}_n)\in \sK^n$ with $\widehat{\bx}_i=\bx_{\sigma(i)}$ and $\widehat{\bz}_i=\bz_{\sigma(i)}$ 
for any $i$. 
Clearly, $\mathrm{id}\in \Pi_{\widehat{\bX}}.$  
Since
$\widehat{\mu}(B_{\mathrm{id}(i)})=\widehat{\mu}(\{i,\ldots, n\})=  \mu(\{\sigma(i),\ldots, \sigma(n)\}) = \mu(B_{\sigma(i)})$
for any $i\in [n]$ and $\widehat{\mu}(B_{\mathrm{id}(n+1)})= \mu(B_{\sigma(n+1)})$, we have 
\begin{align}\label{n:21a}
    \cC_{\mathrm{id}} (\widehat{\bX},\widehat{\mu}) 
    & =  \bigoplus_{i=1}^n \rL(\widehat{\bx}_{i},\widehat{\bx}_{i-1}, \widehat{\mu}(B_{\mathrm{id}(i)}),\widehat{\mu}(B_{\mathrm{id}(i+1)}))\notag
    \\& = \bigoplus_{i=1}^n \rL(\bx_{\sigma(i)},\bx_{\sigma(i-1)},{\mu}({B}_{{\sigma}(i)}),{\mu}({B}_{\sigma(i+1)}))
 =\cC_{\sigma}(\bX,\mu).
\end{align}
Let  $\widehat{\tau}=\sigma^{-1}\circ\tau$ be the product 
of $\sigma^{-1}$ and $\tau$. 
Observe that $\widehat{\tau}\in \Pi_{\widehat{\bZ}}$, as $\widehat{\bz}_{\widehat{\tau}(i)}=\bz_{\tau(i)}$ and $\tau\in\Pi_{\bZ}$. 
Then
\begin{align}\label{n:21b}
\cC_{\widehat{\tau}} (\widehat{\bZ},\widehat{\mu}) & =\bigoplus_{i=1}^n \rL(\widehat{\bz}_{\widehat{\tau}(i)}, \widehat{\bz}_{\widehat{\tau}(i-1)},\widehat{\mu}(B_{\widehat{\tau}(i)}) ,\widehat{\mu}(B_{\widehat{\tau}(i+1)}) )\notag
\\&=\bigoplus_{i=1}^n \rL(\bz_{\tau(i)}, \bz_{\tau(i-1)}, \mu(B_{\tau(i)}), \mu(B_{\tau(i+1)}))=\cC_{\tau}(\bZ,\mu).
\end{align}
By~\ref{ma} and by the well-definedness of the $\cC_{\leqa}$\,-\,operator, we get
$\cC_{\mathrm{id}} (\widehat{\bX},\widehat{\mu}) \leqa \cC_{\widehat{\tau}} (\widehat{\bZ},\widehat{\mu}).$ In consequence, by \eqref{n:21a}\,-\,\eqref{n:21b} follows $\cC_{\sigma}(\bX,\mu)\leqa \cC_{\tau}(\bZ,\mu)$.
To sum up,
we obtain statement \eqref{n:20}  for any $\bX,\bZ\in \sK^n$ 
that differ in only one coordinate and are such that  $\bX\leqA \bZ$.

Put $\bZ^{(i)}=(\bz_1,\ldots, \bz_{i}, \bx_{i+1},\ldots, \bx_n)$ for $i\in [n-1].$
Using $n$-times the above considerations, the monotonicity condition~\ref{m1} 
can be derived as follows 
$$
\cC_{\sigma}(\bX,\mu)\leqa \cC_{\tau^{(1)}}(\bZ^{(1)},\mu)\leqa \ldots \leqa \cC_{\tau^{(n-1)}}(\bZ ^{(n-1)},\mu)\leqa \cC_{\tau}(\bZ,\mu).
$$
Here $\sigma$, $\tau$, and $\tau^{(i)}$ are admissible permutations w.r.t.~$\leqa$ for $\bX$, $\bZ$, and $\bZ^{(i)}$, respectively. 
The proof is complete.
\qed
\end{pf}

\medskip

\begin{pot2}
To make the proof more transparent, we will write the  explicit form of the $\cC_{\leqa}$\,-\,operator,
\begin{align*}
    \cC_{\sigma}(\bX,\mu) = \rL(\bx_{\sigma(1)}, \bzer, 1, \mu(\{\sigma(2)\})) \oplus \rL(\bx_{\sigma(2)}, \bx_{\sigma(1)},  \mu(\{\sigma(2)\}), 0)
\end{align*}
for any $\bX\in \sK^2$ and $\sigma\in \Pi_{\bX}$.

``$\Rightarrow $''
Firstly, we prove (a).
Fix $\bX\in \sK^2$ and $\sigma_1, \sigma_2 \in \Pi_{\bX}$. Clearly, $\bX\leqA[2] \bX$, so  by~\ref{m1} we have
$\cC_{\sigma_1}(\bX, \mu)\leqa \cC_{\sigma_2}(\bX, \mu)$ and  $\cC_{\sigma_2} (\bX, \mu)\leqa \cC_{\sigma_1}(\bX, \mu)$ for any $\mu\in \bM$.
Thus,  due to the antisymmetricity of $\leqa$, we get $\cC_{\sigma_1} (\bX, \mu)=\cC_{\sigma_2}(\bX, \mu)$ for any $\mu$.
By arbitrariness of $\bX,\sigma_1,\sigma_2$, the $\cC_{\leqa}$\,-\,operator is well defined.

We now show (b).
Set $\bX = (\bx, \bv)\in \sK^2$ and $\bZ=(\widehat{\bx}, \bv)\in \sK^2$
such that $\bx \leqa \widehat{\bx} \leqa \bv$.
Clearly, $\bX \leqA[2] \bZ$ and $\mathrm{id} \in\Pi_{\bX}\cap \Pi_{\bZ}$.  It follows from \ref{m1} that 
\begin{align*}
\cC_{\mathrm{id}}(\bX,\mu) & = \rL(\bx, \bzer, 1, \mu(\{2\})) \oplus \rL(\bv, \bx,  \mu(\{2\}),0 ) 
\\&\leqa \rL(\widehat{\bx}, \bzer, 1, \mu(\{2\})) 
    \oplus \rL(\bv, \widehat{\bx},  \mu(\{2\}),0 ) = 
     \cC_{\mathrm{id}}(\bZ,\mu)
\end{align*}
for any $\mu$. By arbitrariness of $\bx \leqa \widehat{\bx} \leqa \bv$, we get point (b).

Put $\bX = (\bu, \bx)\in \sK^2$ and $\bZ=(\bu, \widehat{\bx})\in \sK^2$
such that $\bu \leqa \bx \leqa \widehat{\bx}$.
Clearly, $\bX \leqA[2] \bZ$ and $\mathrm{id} \in\Pi_{\bX}\cap\Pi_{\bZ}$. 
From \ref{m1} we conclude that $\cC_{\mathrm{id}}(\bX,\mu)\leqa \cC_{\mathrm{id}}(\bZ,\mu)$ for any $\mu$.
Lemma~\ref{lem:3}\,(b) 
yields
\begin{align*}
\rL(\bx, \bu, \mu(\{2\}), 0)
\leqa 
\rL(\widehat{\bx}, \bu, \mu(\{2\}), 0)
\end{align*}
for any $\mu$.  By arbitrariness of $\bu \leqa \bx \leqa\widehat{\bx}$, we obtain point (c).

``$\Leftarrow$''
According to Lemma~\ref{lem:10} and the well-definedness of the $\cC_{\leqa}$\,-\,operator, it suffices to show that condition \ref{ma} holds.
Let $\bX=(\bx_1,\bx_2)\in\sK^2$  with $\bx_1\leqa \bx_2$.
Assume that $\bZ = (\bz_1, \bz_2)\in \sK^2$ is such that $\bx_p \lea \bz_p$ for some $p\in [2]$ and $\bx_i = \bz_i$ for $i\neq p$. Clearly, $\bX\leqA[2] \bZ$ and $\Pi_{\bZ}$ contains at most two elements.
We will consider two cases: 
$\bz_1\leqa \bz_2$ and $\bz_2\lea \bz_1$.  Let us start with $\bz_1\leqa \bz_2$.
Evidently, $\mathrm{id} \in \Pi_{\bZ}$.
If $p=1$, we use point (b), and if $p=2$, we use point (c) together with Lemma~\ref{lem:1} to derive the following
\begin{align*}
    \cC_{\mathrm{id}}(\bX,\mu)& = \rL(\bx_{1},\bzer, 1, \mu(\{2\})) \oplus \rL(\bx_{2}, \bx_{1}, \mu(\{2\}), 0)
    \\&\leqa  \rL(\bz_{1}, \bzer, 1, \mu(\{2\})) \oplus \rL(\bz_{2}, \bz_{1}, \mu(\{2\}), 0)= \cC_{\mathrm{id}}(\bZ,\mu),
\end{align*}
so~\ref{ma} is true.   

 We now turn to the case $\bz_2 \lea \bz_1$. 
Then $\Pi_{\bZ}$ consists of a~single element $\tau=(2,1).$  
We show that $p=1$, that is, $\bx_{1}\lea \bz_{1}$ and $\bx_{2}= \bz_{2}$. 
Suppose, contrary to our claim, that 
$p=2$,
i.e.,  $\bx_{1}=\bz_{1}$ and $\bx_{2}\lea  \bz_{2}$. 
Note that $\bx_2 \lea \bz_{2}\lea \bz_{1}=\bx_{1}\leqa \bx_{2}$.  
So by the transitivity of $\lea$ we have $\bx_{2}\lea \bx_{2}$, a~contradiction.  
Therefore $p=1$. 

Define $\bX^{(1)}= (\bx_{2}, \bx_{2})$. 
Clearly,  $\bX \leqA[2] \bX^{(1)}$ and $\Pi_{\bX^{(1)}}=\{\mathrm{id},\tau\}$. By point (b), 
\begin{align}\label{n2:1}
     \cC_{\mathrm{id}}(\bX,\mu)& = \rL(\bx_{1}, \bzer, 1, \mu(\{2\})) \oplus \rL(\bx_{2},\bx_{1}, \mu(\{2\}), 0)\notag
     \\&\leqa \rL(\bx_{2}, \bzer, 1, \mu(\{2\})) \oplus \rL(\bx_{2}, \bx_{2}, \mu(\{2\}), 0)\notag
     \\&= \cC_{\mathrm{id}}(\bX^{(1)},\mu).
\end{align}
As the $\cC_{\leqa}$\,-\,operator is well defined, we conclude that $\cC_{\mathrm{id}}(\bX^{(1)},\mu)=\cC_{\tau}(\bX^{(1)},\mu)$. 
Since $\tau(1)=2$, we have $\bX^{(1)}= (\bz_{\tau(1)}, \bz_{\tau(1)}).$ 
By point (c) and Lemma~\ref{lem:1}, we get 
\begin{align}\label{n2:2}
    \cC_{\tau}(\bX^{(1)},\mu)& = \rL(\bz_{\tau(1)}, \bzer, 1, \mu(\{\tau(2)\})) \oplus \rL(\bz_{\tau(1)}, \bz_{\tau(1)}, \mu(\{\tau(2)\}), 0)\notag
    \\&\leqa \rL(\bz_{\tau(1)}, \bzer, 1, \mu(\{\tau(2)\})) \oplus \rL(\bz_{\tau(2)}, \bz_{\tau(1)}, \mu(\{\tau(2)\}), 0)=\cC_{\tau}(\bZ,\mu).
\end{align}
In consequence, inequalities  \eqref{n2:1}\,-\,\eqref{n2:2} and the 
transitivity of $\leqa$ 
imply $\cC_{\mathrm{id}}(\bX,\mu)\leqa \cC_{\tau}(\bZ,\mu)$, as desired.
\qed
\end{pot2}

\medskip

\begin{pot3} The $\cC_{\leqa}$\,-\,operator  takes the form 
\begin{align*}
    \cC_{\sigma}(\bX,\mu) = 
    \rL(\bx_{\sigma(1)}, \bzer, 1, \mu(B_{\sigma(2)})) 
    \oplus 
    \rL(\bx_{\sigma(2)}, \bx_{\sigma(1)},  \mu(B_{\sigma(2)}), \mu(B_{\sigma(3)}))
    \oplus 
    \rL(\bx_{\sigma(3)}, \bx_{\sigma(2)},  \mu(B_{\sigma(3)}), 0)
\end{align*}
for any $\bX\in \sK^3$ and $\sigma\in \Pi_{\bX}$.

``$\Rightarrow $''
The proof that the $\cC_{\leqa}$\,-\,operator is well defined is analogous to the proof of the same fact in Theorem~\ref{tw:mon_2}, so we omit it.
The validity of statements (b)-(d) can be established by arguments similar to those used in the proof of Theorem~\ref{tw:mon_2}.
To  establish point (b), consider 
$\bX = (\bx, \bv, \bjed)\in \sK^3$ and $\bZ=(\widehat{\bx},\bv, \bjed)\in \sK^3$, where $\bx \leqa \widehat{\bx} \leqa \bv$, and 
$\sigma = \tau =\mathrm{id} \in\Pi_{\bX}\cap \Pi_{\bZ}$. 
To  obtain point (c),  take 
$\bX = (\bu, \bx, \bv)\in \sK^3$ and $\bZ=(\bu, \widehat{\bx}, \bv)\in \sK^3$, where $\bu \leqa \bx \leqa \widehat{\bx}\leqa \bv$, and 
$\sigma = \tau =\mathrm{id} \in\Pi_{\bX}\cap \Pi_{\bZ}$. 
Finally, to  get point (d), set $\bX=(\bzer,\bu, \bx)\in \sK^3$ and $\bZ=(\bzer, \bu, \widehat{\bx})\in \sK^3$, in which $\bu \leqa \bx \leqa \widehat{\bx}$,  and 
$\sigma = \tau = \mathrm{id} \in \Pi_{\bX} \cap \Pi_{\bZ}$.

``$\Leftarrow$''
According to Lemma~\ref{lem:10}, it is sufficient to prove~\ref{ma}, as the $\cC_{\leqa}$\,-\,operator is well defined.
Let $\bX=(\bx_1,\bx_2,\bx_3)\in\sK^3$ with $\bx_1\leqa \bx_2 \leqa \bx_3$.
Assume that $\bZ = (\bz_1, \bz_2, \bz_3)\in \sK^3$ is such that $\bx_p \lea \bz_p$ for some $p\in [3]$ and $\bx_i = \bz_i$ for $i\neq p$. 
Evidently, $\mathrm{id} \in \Pi_{\bX}$ and $\bX \leqA[3] \bZ.$ 
To finish the proof,  it is enough to show that
\begin{align}\label{b:2}
    \cC_{\mathrm{id}}(\mathbf{X}, \mu) \leqa \cC_{\tau}(\mathbf{Z}, \mu) \quad\text{ for some } \tau \in \Pi_{\bZ}.
\end{align}
Consider two cases: 
$\bx_p \lea \bz_p\leqa \bx_{p+1}$ and 
$\bx_{p+1} \lea \bz_p$,  
where, by convention, $\bx_4 = \bjed$. 
 Let us start with the first case.
Put $\tau = \mathrm{id} \in \Pi_{\bZ}$. 
\begin{itemize}[noitemsep]
    \item For $p=1$, using (b) we have 
    \begin{align}
        &\rL(\bx_1, \bzer, 1, \mu(B_{\mathrm{id}(2)})) \oplus \rL(\bx_2, \bx_1,  \mu(B_{\mathrm{id}(2)}), \mu(B_{\mathrm{id}(3)})) \notag \\
        &\quad \leqa \rL(\bz_1, \bzer, 1, \mu(B_{\mathrm{id}(2)})) \oplus \rL(\bx_2, \bz_1,  \mu(B_{\mathrm{id}(2)}), \mu(B_{\mathrm{id}(3)})).\label{n:mon1}
    \end{align}

    \item For $p=2$, employing 
    (c) we obtain 
    \begin{align}
        &\rL(\bx_2, \bx_1, \mu(B_{\mathrm{id}(2)}), \mu(B_{\mathrm{id}(3)})) \oplus \rL(\bx_3, \bx_2,  \mu(B_{\mathrm{id}(3)}), 0)  \notag\\
        &\quad \leqa \rL(\bz_2, \bx_1, \mu(B_{\mathrm{id}(2)}), \mu(B_{\mathrm{id}(3)})) \oplus \rL(\bx_3, \bz_2,  \mu(B_{\mathrm{id}(3)}), 0).\label{n:mon2}
    \end{align}    
    
    \item For $p=3$, applying (d) we get
    \begin{align}\label{n:mon3}
        \rL(\bx_3, \bx_2,  \mu(B_{\mathrm{id}(3)}), 0) \leqa 
        \rL(\bz_3, \bx_2,  \mu(B_{\mathrm{id}(3)}), 0).
    \end{align}
\end{itemize}
Combining \eqref{n:mon1}\,-\,\eqref{n:mon3} with Lemma~\ref{lem:1}  and point (a) 
yields $\cC_{\mathrm{id}}(\bX,\mu) \leqa \cC_{\mathrm{id}}(\bZ, \mu) = \cC_{\tau}(\bZ, \mu).$

Now we examine the case where $\bx_{p+1} \lea \bz_p$.  Then $p \in \{1,2\}$ and $\mathrm{id} \notin \Pi_{\bZ}.$ 
\begin{itemize}[noitemsep]
   \item Let $p=1.$ This means that $\bx_i = \bz_i$ for $i\in \{2,3\}$,  so there are two possibilities: 
   $\bx_2 \lea \bz_1 \leqa \bx_3$ and $\bx_3 \lea \bz_1$. 
   Let us begin with the first. 
   Since $\bz_2 \lea \bz_1 \leqa \bz_3$, $\tau = (2, 1, 3)\in \Pi_{\bZ}$.
    Consider an auxiliary vector $\mathbf{X}^{(1)} = (\bx_2, \bx_2,\bx_3)$. Then $\mathrm{id},\tau \in \Pi_{\mathbf{X}^{(1)} }$ and $\bX \leqA[3] \bX^{(1)}$.
    From point (b), Lemma~\ref{lem:1}, and point (a), we get 
    \begin{align*}
    \cC_{\mathrm{id}}(\bX,\mu)& =
    \rL(\bx_{1}, \bzer, 1, \mu(B_{\mathrm{id}(2)})) 
    \oplus 
    \rL(\bx_{2}, \bx_{1},  \mu(B_{\mathrm{id}(2)}), \mu(B_{\mathrm{id}(3)}))
    \oplus 
    \rL(\bx_{3}, \bx_{2},  \mu(B_{\mathrm{id}(3)}), 0)
    \\&  \leqa
    \rL(\bx_{2}, \bzer, 1, \mu(B_{\mathrm{id}(2)})) 
    \oplus 
    \rL(\bx_{2}, \bx_{2},  \mu(B_{\mathrm{id}(2)}), \mu(B_{\mathrm{id}(3)}))
    \oplus 
    \rL(\bx_{3}, \bx_{2},  \mu(B_{\mathrm{id}(3)}), 0)
    \\&= \cC_{\mathrm{id}}(\mathbf{X}^{(1)}, \mu) = \cC_{\tau}(\mathbf{X}^{(1)}, \mu). 
    \end{align*}
    Since $\bX^{(1)} = (\bz_{\tau(1)}, \bz_{\tau(1)}, \bz_{\tau(3)})$, due to point (c)
    and Lemma~\ref{lem:1}, we obtain
    \begin{align*}
    \cC_{\tau}(\mathbf{X}^{(1)}, \mu)& =  \rL(\bz_{\tau(1)}, \bzer, 1, \mu(B_{\tau(2)})) 
    \oplus 
    \rL(\bz_{\tau(1)}, \bz_{\tau(1)},  \mu(B_{\tau(2)}), \mu(B_{\tau(3)}))
    \oplus 
    \rL(\bz_{\tau(3)}, \bz_{\tau(1)},  \mu(B_{\tau(3)}), 0)
    \\& \leqa 
    \rL(\bz_{\tau(1)}, \bzer, 1, \mu(B_{\tau(2)})) 
    \oplus 
    \rL(\bz_{\tau(2)}, \bz_{\tau(1)},  \mu(B_{\tau(2)}), \mu(B_{\tau(3)}))
    \oplus 
    \rL(\bz_{\tau(3)}, \bz_{\tau(2)},  \mu(B_{\tau(3)}), 0)
     \\& = \cC_{\tau}(\mathbf{Z}, \mu),
    \end{align*}
    which leads us to $\cC_{\mathrm{id}}(\mathbf{X}, \mu) \leqa \cC_{\tau}(\mathbf{Z}, \mu)$, so statement \eqref{b:2} is proved.
    Assume now that $\bx_3 \lea \bz_1$. 
    Since   
    $\bz_2 \leqa \bz_3 \lea \bz_1$,  $\tau = (2, 3, 1)\in \Pi_{\bZ}$. 
    Applying the same reasoning as above results in
    \begin{align*}
        \cC_{\mathrm{id}}(\bX,\mu) &\leqa \cC_{\mathrm{id}}(\bX^{(1)},\mu) = \cC_{\tau^{(1)}}(\bX^{(1)},\mu)\leqa \cC_{\tau^{(1)}}(\bX^{(2)},\mu)
    = \cC_{\tau}(\bX^{(2)},\mu), 
    \end{align*}
    where $\mathbf{X}^{(1)} = (\bx_2,\bx_2,\bx_3)$, $\mathbf{X}^{(2)} =(\bx_3,\bx_2,\bx_3)$, and $\tau^{(1)}=(2, 1, 3)\in \Pi_{\bX^{(1)}}$,
    as $\tau\in \Pi_{\bX^{(2)}}$.
    From point (d) and Lemma~\ref{lem:1} we get $\cC_{\mathrm{id}}(\mathbf{X}, \mu) \leqa \cC_{\tau}(\mathbf{Z}, \mu)$.
    The statement \eqref{b:2} is valid.

    \item Let $p=2$. Then $\bx_3 \lea \bz_2.$ 
    Since $\bz_1 \leqa \bz_3 \lea \bz_2$, 
    $\tau = (1, 3, 2)\in \Pi_{\bZ}$. 
    Take an auxiliary vector $\mathbf{X}^{(1)} = (\bx_1,\bx_3,\bx_3)$ and observe that $\mathrm{id}, \tau\in \Pi_{\mathbf{X}^{(1)}}$, and $\bX \leqA[3] \bX^{(1)}$.  Combining point (c)  with Lemma~\ref{lem:1} we get 
    \begin{align*}
    \cC_{\mathrm{id}}(\bX,\mu) &= \rL(\bx_{1}, \bzer, 1, \mu(B_{\mathrm{id}(2)})) 
    \oplus 
    \rL(\bx_{2}, \bx_{1},  \mu(B_{\mathrm{id}(2)}), \mu(B_{\mathrm{id}(3)}))
    \oplus 
    \rL(\bx_{3}, \bx_{2},  \mu(B_{\mathrm{id}(3)}), 0)
    \\&\leqa 
     \rL(\bx_{1}, \bzer, 1, \mu(B_{\mathrm{id}(2)})) 
    \oplus 
    \rL(\bx_{3}, \bx_{1},  \mu(B_{\mathrm{id}(2)}), \mu(B_{\mathrm{id}(3)}))
    \oplus 
    \rL(\bx_{3}, \bx_{3},  \mu(B_{\mathrm{id}(3)}), 0)
    \\&= \cC_{\mathrm{id}}(\mathbf{X}^{(1)}, \mu).
    \end{align*}
    Since $\bX^{(1)} = (\bz_{\tau(1)}, \bz_{\tau(2)}, \bz_{\tau(2)})$ and $\cC_{\mathrm{id}}(\mathbf{X}^{(1)}, \mu) = \cC_{\tau}(\mathbf{X}^{(1)}, \mu)$ by point (a), we have
    \begin{align*}
    \cC_{\tau}(\mathbf{X}^{(1)}, \mu)& = \rL(\bz_{\tau(1)}, \bzer, 1, \mu(B_{\tau(2)})) 
    \oplus 
    \rL(\bz_{\tau(2)}, \bz_{\tau(1)},  \mu(B_{\tau(2)}), \mu(B_{\tau(3)}))
    \oplus 
    \rL(\bz_{\tau(2)}, \bz_{\tau(2)},  \mu(B_{\tau(3)}), 0).
    \end{align*}
    From point (d) we know that
    $\rL(\bz_{\tau(2)}, \bz_{\tau(2)},  \mu(B_{\tau(3)}), 0)\leqa  \rL(\bz_{\tau(3)}, \bz_{\tau(2)},  \mu(B_{\tau(3)}), 0)$.   
    Using Lemma~\ref{lem:1} we conclude that 
    $
    \cC_{\tau}(\mathbf{X}^{(1)}, \mu)\leqa \cC_{\tau}(\mathbf{Z}, \mu),
    $
    which leads us to $\cC_{\mathrm{id}}(\mathbf{X}, \mu) \leqa \cC_{\tau}(\mathbf{Z}, \mu)$.
    The proof is complete. \qed
\end{itemize}
\end{pot3}

\medskip

\begin{pot4}
``$\Rightarrow$''
The proof that the $\cC_{\leqa}$\,-\,operator is well defined is analogous to the proof of the same fact in Theorem~\ref{tw:mon_2}, so we omit it.
To establish point (b), consider 
$\bX = (\bzer,\ldots,\bzer, \bu,  \bx, \bv, \bjed)\in \sK^n$, $\bZ=(\bzer,\ldots,\bzer, \bu, \widehat{\bx},\bv, \bjed)\in \sK^n$
with $\bu \leqa\bx \leqa \widehat{\bx} \leqa \bv$ and $\sigma=\tau=\mathrm{id}\in \Pi_{\bX}\cap \Pi_{\bZ}$.
To show point (c), set $\bX = (\bzer, \ldots, \bzer, \bu, \bx)\in \sK^n$ and $\bZ=(\bzer,\ldots, \bzer, \bu, \widehat{\bx})\in \sK^n$
with $\bu \leqa \bx \leqa \widehat{\bx}$ and $\sigma=\tau=\mathrm{id}\in \Pi_{\bX}\cap \Pi_{\bZ}$.

``$\Leftarrow$'' 
According to Lemma~\ref{lem:10}, it is sufficient to prove that~\ref{ma} is satisfied.
Let $\bX=(\bx_1,\ldots,\bx_n)\in \sK^n$ be such that $\bx_1\leqa \ldots\leqa \bx_n$ and $\bZ=(z_1,\ldots,z_n)\in \sK^n$ be such that $\bx_p \lea \bz_p$ for some $p\in [n]$ and $\bx_i = \bz_i$ for any $i\neq p$.
Clearly, $\mathrm{id}\in \Pi_{\bX}$ and $\bX \leqA \bZ$. 
Let $q$ be the largest element of the set $\{p,p+1,\ldots,n\}$ such that   $\bx_q\lea \bz_p\leqa \bx_{q+1}$ under the convention $\bx_{n+1}=\bjed$.
We need to prove that
\begin{align}\label{nn:1}   
\cC_{\mathrm{id}}(\bX,\mu) \leqa \cC_{\tau}(\bZ,\mu)\quad \text{ for some } \tau\in \Pi_{\bZ}.
\end{align}

Assume that $q=p$.
In this case we put $\tau = \mathrm{id}\in \Pi_{\bZ}$.
If $1\le p < n$, then by point (b) and Lemma~\ref{lem:1}
\begin{align*}
    \cC_{\mathrm{id}}(\bX,\mu) &= 
    \bigoplus_{i=1}^{p-1} \rL(\bz_i, \bz_{i-1}, \mu(B_{\mathrm{id}(i)}), \mu(B_{\mathrm{id}(i+1)}))
    \oplus
    \rL(\bx_p, \bz_{p-1}, \mu(B_{\mathrm{id}(p)}), \mu(B_{\mathrm{id}(p+1)})) 
    \\&\qquad  \oplus \rL(\bz_{p+1}, \bx_{p}, \mu(B_{\mathrm{id}(p+1)}), \mu(B_{\mathrm{id}(p+2)})) 
    \oplus 
    \bigoplus_{i=p+2}^{n} \rL(\bz_i, \bz_{i-1}, \mu(B_{\mathrm{id}(i)}), \mu(B_{\mathrm{id}(i+1)}))
    \\& \leqa 
    \bigoplus_{i=1}^{p-1} \rL(\bz_i, \bz_{i-1}, \mu(B_{\mathrm{id}(i)}), \mu(B_{\mathrm{id}(i+1)}))
    \oplus
    \rL\big(\bz_p, \bz_{p-1}, \mu(B_{\mathrm{id}(p)}), \mu(B_{\mathrm{id}(p+1)})\big) 
    \\&\qquad  \oplus
    \rL( \bz_{p+1}, \bz_{p}, \mu(B_{\mathrm{id}(p+1)}), \mu(B_{\mathrm{id}(p+2)}) ) 
    \oplus 
    \bigoplus_{i=p+2}^{n} \rL( \bz_i, \bz_{i-1}, \mu(B_{\mathrm{id}(i)}), \mu(B_{\mathrm{id}(i+1)}) )
    \\&=
    \cC_{\mathrm{id}}(\bZ,\mu)=\cC_{\tau}(\bZ,\mu)
\end{align*}  
under the convention ${\tx \bigoplus_{i=1}^0 (\cdot)\oplus \by = \by = \by \oplus \bigoplus_{i=n+1}^{n} (\cdot)}$. 
If $p=n$, then by 
point (c) and Lemma~\ref{lem:1}
\begin{align*}
    \cC_{\mathrm{id}}(\bX,\mu) &= 
    \bigoplus_{i=1}^{n-1} \rL(\bz_i, \bz_{i-1}, \mu(B_{\mathrm{id}(i)}), \mu(B_{\mathrm{id}(i+1)})) \oplus \rL(\bx_n, \bz_{n-1}, \mu(B_{\mathrm{id}(n)}), 0) 
    \\&\leqa \bigoplus_{i=1}^{n-1} \rL(\bz_i, \bz_{i-1}, \mu(B_{\mathrm{id}(i)}), \mu(B_{\mathrm{id}(i+1)})) \oplus \rL(\bz_n, \bz_{n-1}, \mu(B_{\mathrm{id}(n)}), 0)
    \\& = \cC_{\mathrm{id}}(\bZ,\mu)=\cC_{\tau}(\bZ,\mu).
\end{align*}

Assume now that $q\in \{p+1,\ldots, n\}$. 
Since $\bz_i = \bx_i$ for any $i\neq p$ and $\bx_q\lea \bz_p \leqa \bx_{q+1}$, we have
$$
\bz_1 \leqa \ldots \leqa \bz_{p-1} \leqa 
\begin{cases} \bz_{p+1} \leqa \ldots \leqa \bz_q \lea \bz_p \leqa \bz_{q+1} \leqa \ldots \leqa \bz_n & \text{for } q < n, \\
\bz_{p+1} \leqa \ldots  \leqa \bz_q \lea \bz_p & \text{for } q = n.
\end{cases}
$$
Set 
\begin{align*}
    \tau(i) = \begin{cases}
    i   &\text{if } i<p \text{ or } i>q,\\
    i+1 &\text{if } p\le i\le q-1,\\
    p   &\text{if } i=q.
 \end{cases} 
\end{align*}
Evidently, $\tau \in \Pi_{\bZ}$.

Put $\tau^{(0)}=\mathrm{id}$ and $\bX^{(0)}=\bX$. 
To simplify references, we will write $\bx_{p+1}= \bv_1, \ldots, \bx_{q-1}=\bv_{N-1}, \bx_q = \bv_N$, where $N=q-p\ge 1$. 
Clearly, $\bv_1\leqa \ldots\leqa \bv_N$.
We now perform $N$ steps consecutively, in which the $k$-th step, $k\in [N]$, is as follows
\begin{enumerate}[leftmargin=1.4cm]
    \item[\textbf{Step k}] We define 
    $\bX^{(k)} = (\bx_1, \ldots, \bx_{p-1}, \bv_k, \bx_{p+1}, \ldots, \bx_n) \in \sK^n$
    in such a~way that the $p$-th element of $\bX^{(k-1)}$ is replaced with $\bv_k$  and other elements remain the same. 
    Let $\tau^{(k)}$ be defined as follows
    $$
    \tau^{(k)}(i) = 
    \begin{cases}
        i   & \text{if } i<p \text{ or } i>p+k,\\
        i+1 & \text{if } p\le i\le p+k-1,\\
        p   & \text{if } i=p+k.
    \end{cases} 
    $$  
     Then $\tau^{(k-1)}\in \Pi_{\bX^{(k-1)}}\cap \Pi_{\bX^{(k)}}$ and $\tau^{(k)} \in \Pi_{\bX^{(k)}}$. 
    By points (a)-(b) and Lemma~\ref{lem:1},  we get 
    \begin{align}\label{n:3}
        \cC_{\tau^{(k-1)}}(\bX^{(k-1)},\mu)\leqa \cC_{\tau^{(k-1)}}(\bX^{(k)},\mu)=\cC_{\tau^{(k)}}(\bX^{(k)},\mu).
    \end{align}
\end{enumerate}
After $N$ steps, from \eqref{n:3} for $k \in [N]$ we obtain $\cC_{\mathrm{id} }(\bX,\mu)\leqa \cC_{\tau^{(N)}}(\bX^{(N)},\mu)$.
Since $p+N=q$,  $\bX^{(N)}\leqa \bZ$, and $\tau^{(N)}=\tau$, 
from point (c) and Lemma~\ref{lem:1} if $q=n$, or from point (b) and Lemma~\ref{lem:1} if $q<n$, we conclude that 
\begin{align*}
    \cC_{\tau^{(N)}}(\bX^{(N)},\mu) \leqa \cC_{\tau}(\bZ,\mu).
\end{align*}
Therefore statement \eqref{nn:1} is valid. \qed
\end{pot4}

\section*{Appendix C}

Let us recall the construction of the IV dissimilarity function proposed in \cite[Prop.~4.8]{takac2022}. To do this, we will first  introduce the basic notions.

For $\bx = [x^l, x^u] \in \sLi$, define $K_{\alpha}(\bx) = (1-\alpha) x^l + \alpha x^u$ for any $\alpha \in [0,1]$, $w(\bx) = x^u - x^l$,  
and
\begin{align}\label{dis:0}
    \lambda_{\alpha} (\bx) = \frac{w(\bx)}{\frac{K_{\alpha}(\bx)}{\alpha} \wedge \frac{1-K_{\alpha}(\bx)}{1-\alpha}}, \quad  \alpha\in (0,1)
\end{align}
under convention $\tfrac{0}{0}=0$ \cite[Eqs.~(6) and (7)]{takac2022}. Hereafter, $a\wedge b=\min\{a,b\}$ for any $a,b\in\mR_+$.

\begin{proposition}\cite[Prop.~4.18]{takac2022}\label{prop:takac}
Let $\alpha,\beta \in (0,1)$ be such that $\alpha \neq \beta$. 
 Let $M_{\rd}\colon [0,1]^2 \to [0,1]$ be a~symmetric  aggregation function, $\delta_{\rd}\colon [0,1]^2 \to [0,1]$ be a~strictly monotone  dissimilarity function. Then the function $\rd\colon (\sLi)^2 \to \sLi$ defined by
\begin{align*}
    \rd(\bx, \by) = \bz , \qquad \text{ where } \qquad \begin{cases}
        K_{\alpha}(\bz) = \delta_{\rd}(K_{\alpha}(\bx), K_{\alpha}(\by)),\\
        \lambda_{\alpha}(\bz)  = M_{\rd}(\lambda_{\alpha}(\bx), \lambda_{\alpha}(\by))
    \end{cases}
\end{align*}
is IV dissimilarity function w.r.t.~$\leqab$. Moreover,  
\begin{align}\label{tak:1}
    \bz = \big[K_{\alpha}(\bz) - \alpha w(\bz),\; K_{\alpha}(\bz) + (1-\alpha)w(\bz)],
\end{align}
where $w(\bz)$ is given by  \eqref{dis:0} with~$\bx=\bz$.
\end{proposition}

We will demonstrate that the IV dissimilarity function constructed according to Proposition~\ref{prop:takac} does not satisfy point (b)  of  Corollary~\ref{cor:4.14} 
 provided  $\delta_{\rd}([0,1],0) = [0,1]$.

Let $\alpha$, $\beta$, $M_{\rd}$, and $\delta_{\rd}$ be such as in Proposition~\ref{prop:takac}.
Assume additionally that $\delta_{\rd}([0,1], 0) = [0,1]$.
Put $\bx_1 = [0, x_1]$ and $\bx_2 = [0, x_2]$ with $x_1,x_2 \in (0,1]$ and $x_1 < x_2$. 
Hence $\bx_1 \leqab \bx_2$ for any $x_1<x_2$.
Let $\bz_1 = \rd(\bx_1, \bzer)$, $\bz_2 = \rd(\bx_2, \bzer)$, and $\bz_{12} = \rd(\bx_1, \bx_2)$ take the form \eqref{tak:1}.
Suppose that point (b) of Corollary~\ref{cor:4.14} is true, i.e.,
\begin{align}\label{tak:2}
    \bz_1 \oplus_{\mathrm{iv}} \bz_{12} = \bz_2
\end{align}
for any $x_1 <x_2$.
We establish  $K_{\alpha}(\bz_j)$ and $\lambda_{\alpha}(\bz_j)$ for all $j\in \{1, 2, 12\}$: 

\begin{itemize}[noitemsep]
    \item $K_{\alpha}(\bz_i)= \delta_{\rd}(K_{\alpha}(\bx_i), 0) = \delta_{\rd}(\alpha x_i, 0)$ for $i\in [2]$,
    
    \item $K_{\alpha}(\bz_{12})= \delta_{\rd}(K_{\alpha}(\bx_1), K_{\alpha}(\bx_2)) = \delta_{\rd}(\alpha x_1, \alpha x_2)$,

    \item $\lambda_{\alpha}(\bz_i)= M_{\rd}(\lambda_{\alpha}(\bx_i), \lambda_{\alpha}(\bzer)) = M_{\rd}(1,0)$ for $i\in [2]$,
    
    \item $\lambda_{\alpha}(\bz_{12})= M_{\rd}(\lambda_{\alpha}(\bx_1), \lambda_{\alpha}(\bx_2)) = M_{\rd}(1,1)=1$.
\end{itemize}
Put $A_i = \delta_{\rd}(\alpha x_i, 0) $ for $i\in [2]$ and $A_{12}=  \delta_{\rd}(\alpha x_1, \alpha x_2)$.
Thus
\begin{align*}
    \bz_j &= \big[ A_j - \alpha w(\bz_j), \; A_j + (1-\alpha) w(\bz_j)\big], \quad j \in \{1, 2, 12\}.
\end{align*}
By equating the left and right ends of the intervals in equation \eqref{tak:2},  
we obtain
\begin{align*}
    \begin{cases}
    \alpha (w(\bz_1)+w(\bz_{12}) - w(\bz_2))  = A_1 + A_{12} - A_2,\\
    (1-\alpha) (w(\bz_1) + w(\bz_{12}) - w(\bz_2)) = - (A_1 + A_{12} - A_2).
    \end{cases}
\end{align*}
In consequence, we have $A_2 = A_1 + A_{12}$ and
\begin{align}\label{tak:4}
    w(\bz_2) = w(\bz_1) + w(\bz_{12})
\end{align}
for any $x_1 < x_2$.
Let us now determine the values of $w(\bz_i)$ using \eqref{dis:0}:
\begin{itemize}[noitemsep]
    \item ${\tx w(\bz_i)= \lambda_{\alpha}(\bz_i) \cdot (\frac{K_{\alpha}(\bz_i)}{\alpha} \wedge \frac{1-K_{\alpha}(\bz_i)}{1-\alpha}) = M_d(1,0) ( \frac{A_i}{\alpha} \wedge \frac{1-A_i}{1-\alpha})}$ for $i\in [2]$,

    \item ${\tx w(\bz_{12})= \frac{A_{12}}{\alpha} \wedge \frac{1-A_{12}}{1-\alpha}}$.
\end{itemize}
Since $A_{12} = A_2 - A_1 >0$, as $\delta_{\rd}$ is strictly monotone, \eqref{tak:4} takes the form
\begin{align}\label{tak:3}
    M_{\rd}(1,0) \Big(\frac{A_2}{\alpha} \wedge \frac{1- A_2}{1-\alpha}\Big) = M_{\rd}(1,0) \Big(\frac{A_1}{\alpha} \wedge \frac{1- A_1}{1-\alpha}\Big)  + \Big( \frac{A_2- A_1}{\alpha} \wedge \frac{1- A_2 + A_1}{1-\alpha}\Big)
\end{align}
for any $x_1 < x_2$.
Note that
\begin{align}\label{tak:5}
\begin{cases}
        \frac{A}{\alpha} \le \frac{1- A}{1-\alpha} \quad \Leftrightarrow\quad A \le \alpha,\\
    \frac{A_2- A_1}{\alpha} \le  \frac{1- A_2 + A_1}{1-\alpha}\quad \Leftrightarrow\quad A_2 - A_1 \le \alpha. 
\end{cases}
\end{align}

Given that $\delta_{\rd}([0,1],0)=[0,1]$, there exist $x_1$ and $x_2$ such that
$A_2 \le \alpha$ or $A_1 \le \alpha < A_2.$ 
Assume
that $A_2 \le \alpha$ for some $x_1 < x_2$. 
Then $A_1 \le \alpha$ and $A_2 - A_1 \le \alpha$. 
By~\eqref{tak:5}, \eqref{tak:3} has the form
\begin{align*}
    M_{\rd}(1,0) \frac{A_2}{\alpha} = M_{\rd}(1,0) \frac{A_1}{\alpha} + \frac{A_2- A_1}{\alpha}.
\end{align*}
Hence $M_{\rd}(1,0)=1$.
Assume that $A_1 \le \alpha < A_2$ for some $x_1 < x_2$. 
If $A_2 - A_1 \le \alpha$, then \eqref{tak:3} can be rewriten as
\begin{align*}
    \frac{1-A_2}{1-\alpha} =  \frac{A_1}{\alpha} + \frac{A_2-A_1}{\alpha},
\end{align*}
so $A_2 = \alpha$,
as $M_{\rd}(1,0)=1$, a~contradiction.
If $A_2 - A_1 > \alpha$, then \eqref{tak:3} is  of the form
\begin{align*}
    \frac{1-A_2}{1-\alpha} =  \frac{A_1}{\alpha} +  \frac{1- A_2 + A_1}{1-\alpha} 
\end{align*}
and implies $A_1 = 0$,
which contradicts $A_1>0$, as $x_1 >0$.

        

        


\end{document}